\documentclass[envcountsame,envcountsect]{svmult}

\usepackage{makeidx}         
\usepackage{graphicx}        
\usepackage{multicol}        
\usepackage[bottom]{footmisc}
\usepackage{amsmath,amssymb,latexsym} 

\makeindex
\smartqed
\spnewtheorem{construction}[theorem]{Construction}{\bfseries}{\rmfamily}

\DeclareMathOperator\ad{ad}

\DeclareMathOperator\Aut{Aut}
\DeclareMathOperator\CDer{CDer}
\DeclareMathOperator\Ct{Cent}

\DeclareMathOperator\Der{Der}
\DeclareMathOperator\diag{diag}
\DeclareMathOperator\ev{ev}
\DeclareMathOperator\Ext{Ext}
\DeclareMathOperator\GL{GL}

\DeclareMathOperator\End{End}
\DeclareMathOperator\Hom{Hom}
\DeclareMathOperator\id{id}
\DeclareMathOperator\Id{Id}

\DeclareMathOperator\LT{LT} 
\DeclareMathOperator\Mat{Mat} 

\DeclareMathOperator\pe{pe}
\DeclareMathOperator\rad{rad}
\DeclareMathOperator\rank{rank}
\DeclareMathOperator\SCDer{SCDer}
\DeclareMathOperator\spann{span}
\DeclareMathOperator\speciallinear{sl}
\DeclareMathOperator\supp{supp}
\DeclareMathOperator\ssp{ssp}
\DeclareMathOperator\trace{tr}
\DeclareMathOperator\TKK{TKK}

\newcommand\andd{\quad \text{and} \quad}

\newcommand\form{(\,\, \vert\,\, )}

\newcommand\order[1]{\left\vert #1 \right\vert}
\newcommand\ot{\otimes}
\newcommand \set[1]{\{ \, #1 \, \}}
\newcommand\suchthat{\mid}

\newcommand \fg {\mathfrak g}
\newcommand \fh {\mathfrak h}

\newcommand \fs {\mathfrak s}

\newcommand\rA{A}

\newcommand\rC{{C}}
\newcommand\rD{{D}}
\newcommand\rE{{E}}
\newcommand\rH{{H}}

\newcommand\rL{{L}}
\newcommand\rM{{M}}

\newcommand\rR{{R}}

\newcommand\rT{T}
\newcommand\rV{V}
\newcommand\rX{{X}}

\newcommand\cC{{\mathcal C}}
\newcommand\cD{{\mathcal D}}
\newcommand\cP{{\mathcal P}}

\newcommand{\al}{\alpha}
\newcommand{\ep}{\varepsilon}
\newcommand{\gm}{\gamma}
\newcommand{\Gm}{\Gamma}
\newcommand{\lm}{\lambda}
\newcommand{\Lm}{\Lambda}
\newcommand{\ph}{\varphi}
\newcommand{\sg}{\sigma}

\newcommand\bq{{\mathbf q}}
\newcommand\be{{\mathbf e}}
\newcommand\bsg{{\boldsymbol{\sg}}}

\newcommand{\blm}{{\bar\lambda}}
\newcommand{\bLm}{{\bar\Lambda}}
\newcommand\bmu{{\bar \mu}}

\newcommand\bbC{\mathbb C}
\newcommand\bbZ{\mathbb Z}
\newcommand\Zn{{\mathbb Z^n}}

\newcommand{\De}{\Delta}
\newcommand{\Dec}{{\Delta^\times}}

\newcommand{\Dei}{{\Delta_{\textrm{ind}}}}

\newcommand\base{{\mathbb F}} 
\newcommand\bi{\simeq_{bi}}
\newcommand\coreE{\rE_\text{c}}
\newcommand\ccoreE{\rE_\text{cc}}
\newcommand\eu{e} 
\newcommand\grd{{\text{gr}*}} 
\newcommand\grtwist{\text{gr}}
\newcommand\hi{{(h)}}
\newcommand\ig{\simeq_{ig}}
\newcommand\opalg{\text{op}}
\newcommand\shom{s}
\newcommand\rLs{\rL^{(\shom)}}

\newcommand\si{{(\shom)}}
\newcommand\sip{{(\shom')}}

\newcommand\ui{{(u)}}
\newcommand\Ztwo{\bbZ_2}

\newcommand\ass[1]{(\ref{con:famEALA}#1)}

\begin{document}

\title*{Isotopy for extended affine Lie algebras and Lie tori}
\titlerunning{Isotopy for EALAs and Lie tori}
\author{Bruce Allison\inst{1}\thanks{Bruce Allison
gratefully acknowledges the support
of the Natural Sciences and Engineering Research Council of Canada.}\and
John Faulkner\inst{2}}
\institute{Department of Mathematics and Statistics, University of
Victoria, Victoria BC V8W~3P4 Canada
\texttt{ballison@uvic.ca}
\and Department of Mathematics,
University of Virginia, Charlottesville VA 22904-4137 USA
\texttt{jrf@virginia.edu}}


\maketitle

\abstract{Centreless Lie tori have been used by E.~Neher to construct all extended affine Lie algebras (EALAs).
In this article, we study isotopy for centreless Lie tori, and show that Neher's construction provides a
1-1 correspondence between centreless Lie tori up to isotopy and families of EALAs up to isomorphism.
Also, centreless Lie tori can be coordinatized by unital algebras that are in general nonassociative,
and, for many types of centreless Lie tori, there are classical definitions of isotopy
for the coordinate algebras. We show for those types that an isotope of a Lie torus is
coordinatized by an isotope of its coordinate algebra, thereby connecting the two notions of isotopy.
In writing the article, we have not assumed prior knowledge of the theories of EALAs, Lie tori or isotopy.
}


\section{Introduction}
\label{sec:introduction}

If $\rA$  is a  unital associative algebra and $u$ is invertible in $\rA$,
one can define an algebra $\rA^\ui$, called the $u$-isotope of $\rA$,
which is equal to $\rA$
as a vector space but has a new product
$x\cdot_u y = xuy$.
This  isotope $\rA^\ui$ is again unital and associative but with a shifted
identity element $u^{-1}$.  More generally there are definitions of isotope
for several other classes of unital nonassociative algebras, notably
Jordan algebras  \cite[\S I.3.2]{Mc2}, alternative algebras \cite{Mc1} and
associative algebras with involution
\cite[\S I.3.4]{Mc2}.
In each case, the isotope is obtained very roughly by shifting the identity element in the algebra,
and two algebras are said to be  isotopic if one is isomorphic to an isotope of the other.
In the associative case, the $u$-isotope $\rA^\ui$ is
isomorphic to $\rA$ under left multiplication by $u$, and therefore
isotopy has not played a role in associative theory.
That is not true in general though, and in particular
isotopy plays an important role in Jordan theory (see for example \cite[\S II.7]{Mc2}).

In contrast, isotopes and isotopy have not
been defined for Lie algebras, for the evident reason that Lie algebras are not
unital.\footnote{There is a notion of isotopy that makes sense for nonunital algebras and
hence Lie algebras \cite{A}, but it is quite different from the one discussed in this article.}
In this article we study notions of isotope and isotopy,
which were recently introduced in \cite{ABFP2}  for a class of graded Lie algebras called Lie
tori.
The point to emphasize here is that forming an isotope of a Lie torus
does not change the multiplication at all, but rather it shifts the grading.

We are primarily interested in the case when the Lie torus is centreless,
and there are two basic reasons why we are interested in isotopy in this case.
First, centreless Lie tori arise naturally in the construction of families
of extended affine Lie algebras (EALAs), and we see
in this article that isotopes and isotopy play a natural and fundamental role in the theory
of both of these classes of Lie algebras.
In fact, we show in section \ref{sec:isotopyEALA} using
a construction of E.~Neher, that there is a 1-1 correspondence
between centreless Lie tori up to isotopy and families of EALAs up to bijective
isomorphism.  (This terminology is explained in section \ref{sec:isotopyEALA}.)

Second, any Lie torus is by definition graded by the group $Q\times \Lm$, where
$Q$ is the root lattice of a finite irreducible root system $\De$ and $\Lm$ is
(in this article) a finitely generated free abelian group.  Depending on the type of $\De$,
centreless Lie tori are coordinatized by unital algebras that are in general nonassociative.
For example, in types A$_1$, A$_2$, A$_r\ (r\ge 3)$
and C$_r\ (r\ge 4)$ the coordinate algebras are respectively
Jordan algebras, alternative algebras, associative algebras and
associative algebras with involution.
It turns out that an isotope of a centreless Lie torus is
coordinatized by an isotope of its coordinate algebra,
thereby connecting our notion of isotopy for Lie tori
with the classical notions for unital nonassociative algebras.
This fact can then be used to obtain necessary and sufficient conditions
for two centreless Lie tori of the same type to be isotopic
(in terms of their coordinate algebras).
We describe these results in detail for types A$_1$, A$_2$, A$_r\ (r\ge 3)$
and C$_r\ (r\ge 4)$ in sections \ref{sec:A1}, \ref{sec:A2}, \ref{sec:Ar},
and \ref{sec:Cr}, and we remark very briefly on the other types
in section \ref{sec:conclusions}.

Section  \ref{sec:conclusions} also includes a brief discussion of an alternate
approach, using multiloops, to constructing centreless Lie tori and studying isotopy.

In writing this article, we have not assumed
prior knowledge of the theories of EALAs, Lie tori or isotopy.  In fact, we
hope that this article
will serve as a useful introduction to these theories and their interconnections.
For this reason, we begin in sections \ref{sec:EALA}, \ref{sec:Lietori}, \ref{sec:construction}
and \ref{sec:isotopy} by recalling some of the important definitions and facts about
EALAs, Lie tori and isotopy for Lie tori.
For the same reason, we have included a brief discussion
of coordinatization in section \ref{sec:coordinates}; and
in sections \ref{sec:A1} and \ref{sec:A2} we have included
some arguments that are familiar to experts in Jordan and alternative algebra
theory, but may be less familiar to others.

\smallskip\noindent \emph{Acknowledgements:}
We thank  Stephen Berman, Erhard Neher and Arturo Pianzola
for stimulating conversations about the topics of this paper.  We also
thank Erhard Neher for his helpful comments about a preliminary version of this article.

\smallskip\noindent
\emph{Assumptions and notation:}
Throughout this work \emph{we assume that $\base$ is a field of
characteristic $0$}. All algebras are assumed to be algebras over
$\base$, and  all algebras (except Lie algebras) are assumed
to be unital.  We also assume that \emph{$\Lm$ is a finitely generated
free abelian group}.  We denote  the rank of $\Lm$, which is an integer $\ge 0$,
by $\rank(\Lm)$.

If $\rE$ is a Lie algebra and $\rH$ is an ad-diagonalizable
subalgebra of $\rE$, we let
$\rE_\rho = \set{x\in \rE \suchthat [h,x] = \rho(h)x \text{ for } h\in \rH}$ for
$\rho$ in the dual space $\rH^*$ of $\rH$.
Then, we define
$\De(\rE,\rH) = \set{\rho\in \rH^* \suchthat \rE_\rho \ne 0}$\index{De(E,H)@$\De(\rE,\rH)$}
and we call elements of  $\De(\rE,\rH)$ \emph{roots of $\rE$ relative to $\rH$};
so in particular, $0$ is a root (if $\rE \ne 0$).

\smallskip\noindent
\emph{Some terminology for graded algebras:}
If $\rA = \oplus_{\lm\in\Lm} \rA^\lm$ is a $\Lm$-graded algebra, the $\Lm$-\emph{support} of $\rA$
is $\supp_\Lm(\rA) = \set{\lm\in\Lm\suchthat \rA^\lm\ne 0}$.
If $\Gm$ is a subgroup of $\Lm$ and $\rA$ is a $\Gm$-graded algebra or vector space,
we regard $\rA$ as $\Lm$-graded by setting $\rA^\lm = 0$ for  $\lm\in\Lm\setminus\Gm$.
If $\rA$ is a $\Lm$-graded algebra and $\rA'$ is a
$\Lm'$-graded algebra, we say that $\rA$ and $\rA'$ are
\emph{isograded-isomorphic}\index{isograded-isomorphic}
\index{isograded-isomorphic!$\eta_\textrm{gr}$},
written $\rA\ig \rA'$,  if there exists an algebra isomorphism $\eta: \rA \to \rA'$ and a group
isomorphism
$\eta_\textrm{gr} : \Lm \to \Lm'$
such that
$\eta(\rA^\lm) = \rA'^{\eta_\grtwist(\lm)}$ for $\lm\in \Lm$; in  that case
$\eta_\grtwist$ is uniquely determined by $\eta$ if
$\langle\supp_\Lm(\rA)\rangle= \Lm$.
(Here $\langle\supp_\Lm(\rA)\rangle$ denotes the subgroup of $\Lm$
generated by $\supp_\Lm(\rA)$.)
If $\rA$ and $\rA'$ are
$\Lm$-graded algebras, we say that $\rA$ and $\rA'$ are
\emph{graded-isomorphic}\index{graded-isomorphic},
written $\rA \simeq_\Lm \rA'$,
if there is an algebra isomorphism $\eta: \rA\to \rA'$ that
preserves the grading (that is $\eta_\grtwist = \id$).
A \emph{$\Lm$-graded algebra with involution} is a pair
$(\rA,\iota)$, where $\rA$ is a $\Lm$-graded algebra and $\iota$
is an involution (antiautomorphism of period 2) of $\rA$ that is graded.
There is an evident extension of the terms isograded-isomorphic and graded-isomorphic
and the notations $\ig$ and $\simeq_\Lm$
for graded algebras with involution.  Finally, the \emph{group algebra} $\base[\Lm] = \oplus_{\lm\in \Lm} \base t^\lm$
of $\Lm$  with $t^\lm t^{\lm'} =  t^{\lm + \lm'}$ is naturally a $\Lm$-graded algebra.  If we choose a basis
$\set{\lm_1,\dots,\lm_n}$ for $\Lm$, then $\base[\Lm]$ is isograded-isomorphic to the
algebra $\base[t_1^{\pm 1},\dots,t_n^{\pm 1}]$ of \emph{Laurent polynomials} with its   natural
$\Zn$-grading.

\section{Extended affine Lie algebras}
\label{sec:EALA}

The definition of an EALA has evolved from \cite{HT},
where these algebras were introduced under the name
elliptic irreducible quasi-simple Lie algebra, and from \cite{BGK} and \cite{AABGP}, where many
of their properties were developed.  We will use the definition
given in \cite{N2}, which has the important advantage that it makes sense
over any field of characteristic~0.  Some of the facts about EALAs mentioned
below were proved in the setting of \cite{AABGP}, but they can be verified in a similar
fashion in our more general setting.

Recall  that an
\emph{extended affine Lie algebra}\index{extended affine Lie algebra}
\index{EALA|see{extended affine Lie algebra}}
(EALA) is a triple
\linebreak $(\rE,\form,\rH)$, where $\rE$ is a Lie algebra over $\base$,
$\form$ is a nondegenerate invariant symmetric
bilinear form on $\rE$, and $\rH$ is a finite dimensional nonzero
self-centralizing ad-diagonalizable subalgebra of $\rE$, such that
a list of axioms labeled as (EA3)--(EA6) are satisfied
\cite{N2}.  For the purposes of this  article we do not need the precise
statement of these axioms, except to say that they are modeled after the properties of
finite dimensional split simple Lie algebras and
affine Kac-Moody Lie algebras.
If $(\rE,\form,\rH)$ is an EALA, we also say that \emph{$\rE$ is an EALA with respect to $\form$ and $H$},
or simply  that \emph{$\rE$ is an EALA}. Roots of $\rE$ relative to $\rH$
are simply called \emph{roots} of $\rE$.

If $\rE$ is an EALA, we can
as usual  transfer the restriction of
$\form$ to $\rH$ to a nondegenerate
form $\form$ on the dual space $\rH^*$.
If $R = \De(\rE,\rH)$, $\rV = \spann_\base(R)$, and
$\rad(\rV)$ is the radical of the restriction of $\form$ to $\rV$,
then the image of $R$ in $\rV/\rad(\rV)$ is a finite irreducible
root system (including 0) whose
type is called the
\emph{type}\index{extended affine Lie algebra!type}
 of $\rE$. (See the beginning of \S \ref{sec:Lietori} below for our conventions
about finite irreducible root systems.)

A root $\rho$ of an EALA $\rE$ is called
\emph{isotropic}\index{extended affine Lie algebra!isotropic root}
if $(\rho,\rho) = 0$,
and otherwise called \emph{nonisotropic}.
According to one of the axioms for an EALA, the group generated by the isotropic roots
is a free abelian group of finite rank, and its rank is called the
\emph{nullity}\index{extended affine Lie algebra!nullity}
of $\rE$.
Then the extended affine Lie algebras of nullity 0 and 1 are precisely the
finite dimensional split simple Lie algebras and the
affine Kac-Moody Lie algebras respectively. (See \cite[Remark 1.2.4]{ABFP2} for nullity 0
and \cite{ABGP} for  nullity 1.)

If $(\rE,\form,\rH)$ is an EALA, then so  is $(\rE,a\form,\rH)$ for $a\in \base^\times$,
and it is sometimes convenient to adjust the form in this way.
For this reason, the following is a natural notion of isomorphism for EALAs.

\begin{definition}
\label{def:EALAisom}
If  $(\rE,\form,\rH)$ and $(\rE',\form',\rH')$ are  EALAs,
an
\emph{isomorphism}\index{extended affine Lie algebra!isomorphic}
from $(\rE,\form,\rH)$ onto $(\rE',\form',\rH')$
is a Lie algebra isomorphism $\chi : \rE \to \rE'$ such that $(\chi(x)\vert\chi(y))' = a(x\vert y)$
for some $a\in \base^\times$ and all $x,y\in \rE$, and $\chi(\rH) = \rH'$.
If such a map exists we say that $(\rE,\form,\rH)$   and $(\rE',\form',\rH')$ are \emph{isomorphic}.
For short we often say that $\chi$ is an \emph{EALA isomorphism} from $\rE$ onto $\rE'$
and that $\rE$ and $\rE'$ are \emph{isomorphic as EALAs}.
\end{definition}

\begin{definition}
\label{def:core}
If $\rE$ is an EALA,
the
\emph{core}\index{extended affine Lie algebra!core}
of $\rE$ is the
subalgebra $\coreE$ of $\rE$ that is generated by the root spaces of
$\rE$ corresponding to nonisotropic roots, and the
\emph{centreless core}\index{extended affine Lie algebra!centreless core}
of $\rE$  is
$\ccoreE := \coreE/Z(\coreE)$.\footnote{We denote the centre
of a Lie algebra $\rL$ by $Z(\rL)$.  $\rL$ is called
\emph{centreless} if $Z(\rL) = 0$.}
\end{definition}

\section{Lie tori}
\label{sec:Lietori}

With the construction of families of EALAs in mind, Yoshii gave a definition of a
Lie torus in \cite{Y6}.  An equivalent definition, which we recall next,
was given by Neher in~\cite{N1}.

It will  be convenient
for us to work with root systems that contain 0.
So by a \emph{finite irreducible root system} we will mean
a finite subset $\De$ of a finite dimensional vector space
$\rX$ over $k$ such that $0\in \De$ and
$\Dec := \De \setminus\set{0}$ is
a finite irreducible root system in $\rX$ in the usual sense
(see \cite[Chap VI, \S 1, no.~1, Definition~1]{B2}).
With this convention we will \emph{assume for the rest of this section
that   $\De$ is a finite irreducible root system in
a finite dimensional vector space $\rX$ over $k$}.
Recall that $\De$ is  said to be
\emph{reduced}\index{reduced root system}
if $2 \al \notin \Dec$ for $\al \in \Dec$.
If $\De$ is reduced then
$\De$ has type A$_\ell (\ell \ge 1)$, B$_\ell (\ell \ge 2)$,
C$_\ell (\ell \ge 3)$, D$_\ell (\ell \ge 4)$, E$_6$, E$_7$,
E$_8$, F$_4$, or G$_2$, whereas if $\De$ is not reduced,
$\De$ has type BC$_\ell (\ell \ge 1)$
\cite[Chapter VI, \S 4 and 14]{B2}.
We let
\[\Dei := \set{0}\cup \set{\al\in \Dec \suchthat
\textstyle \frac 12 \al\notin \De},\]
in which case $\Dei$ is an irreducible reduced root system in
$\rX$, and $\Dei = \De$ if $\De$ is reduced.
Also we let
\[Q = Q(\De) := \spann_\bbZ(\De)\]
be the   \emph{root lattice} of $\De$.

If $\rL = \textstyle\bigoplus_{(\al,\lm)\in Q\times\Lm}\rL_{\al}^{\lm}$
is a $Q\times\Lm$-graded algebra, then
$\rL = \bigoplus_{\lm \in \Lm} \rL^\lm$ is
$\Lm$-graded and $\rL = \bigoplus_{\al \in Q} \rL_\al$ is $Q$-graded,
with
\[\rL^\lm =\textstyle \bigoplus_{\al\in Q} \rL_{\al}^{\lm} \text{ for } \lm\in \Lm
\quad \andd \quad \rL_\al = \bigoplus_{\lm\in \Lm} \rL_\al^\lm \text{
for }\al\in Q.\]

\begin{definition}
\label{def:Lietorus}
A
\emph{Lie $\Lm$-torus of type $\Delta$}\index{Lie torus}
is a $Q\times\Lm$-graded
Lie algebra $\rL = \bigoplus_{(\al,\lm)\in Q\times\Lm}\rL_{\al}^{\lm}$
over $k$ which satisfies:
\begin{description}[(LT4) ]
\item[(LT1)]  $\rL_\al = 0$ for $\al\in Q\setminus \De$.
\item[(LT2)]
\begin{description}[(ii) ]
\item[(i)] If $0\ne \al\in \De$, then $\rL_\al^0 \ne 0$.
\item[(ii)]
If $0\ne \al\in Q$, $\lm\in\Lm$ and $\rL_\al^\lm\ne 0$,
then there exist elements $e_\al^\lm\in \rL_\al^\lm$
and $f_\al^\lm\in \rL_{-\al}^{-\lm}$ such that
\[\rL_\al^\lm = \base e_\al^\lm,\quad \rL_{-\al}^{-\lm} = \base f_\al^\lm,\]
and
\begin{equation}
\label{eq:basic}
[[e_\al^\lm,f_\al^\lm],x_\beta] = \langle\beta,\al^\vee\rangle x_\beta
\end{equation}
for $x_\beta\in \rL_\beta$,  $\beta\in Q$.
\end{description}
\item[(LT3)] $\rL$ is generated as an algebra by the spaces $\rL_\al$, $\al\in \Dec$.
\item[(LT4)] $\Lm$ is generated as a group by $\supp_\Lm(\rL)$.
\end{description}
In that case the $Q$-grading of $\rL$ is called
the
\emph{root grading}\index{Lie torus!root grading and external grading}
of $\rL$, and the $\Lm$-grading
of $\rL$ is called the
\emph{external grading}
of $\rL$.
If $\De$ has type X$_\ell$, we also say  that
\emph{$\rL$ has type~X$_\ell$}\index{Lie torus!type}.
\end{definition}

\begin{remark}
\label{rem:LT5}  Suppose $\rL$ is a Lie $\Lm$-torus of  type $\De$.

(i) It is shown in \cite[Prop.~1.1.10]{ABFP2} that
$\supp_Q(\rL)$ equals either  $\De$ or~$\Dei$.

(ii) It is sometimes convenient to assume the following additional axiom:
\begin{description}[(LT5) ]
\item[(LT5)]  $\supp_Q(\rL) = \De$.
\end{description}
Note that by (i), (LT5) holds automatically if $\De$ is reduced.
Also, if (LT5)  does not hold, then we can
replace $\De$ by $\Dei$, in which case (LT5) holds.
Thus, when convenient there is no loss of generality in assuming  (LT5).
\end{remark}

In the study of  Lie tori it is not convenient to fix a particular realization of the
root system $\De$ or a particular identification of the group $\Lm$ with $\Zn$.
For this reason, the following is a natural notion of isomorphism for Lie tori.

\begin{definition}\label{def:bi-isomorphism}\
If $\rL$ is  a Lie  $\Lm$-torus of type $\De$  with $Q = Q(\De)$
and $\rL'$ is a Lie $\Lm'$-torus of type $\De'$ with $Q'= Q(\De')$,
$\rL$ and $\rL'$ are said to be
\emph{bi-isograded-isomorphic}\index{Lie torus!bi-isograded-isomorphic}
\index{Lie torus!bi-isograded-isomorphic!$\ph_r$ and $\ph_e$}
or \emph{bi-isomorphic}\index{Lie torus!bi-isomorphic}
for short,  if there is an algebra isomorphism $\ph : \rL \to \rL'$,
a group isomorphism
$\ph_r: Q \to Q$,
and a group isomorphism $\ph_e : \Lm \to \Lm'$
such that
\[\ph(\rL_\al^\lm) = {\rL'}_{\ph_r(\al)}^{\ph_e(\lm)}\]
for $\al\in Q(\De)$ and $\lm\in \Lm$.  In that case
we write
\[\rL\bi\rL'\]
and we call $\ph$ a \emph{bi-isograded-isomorphism},
or a \emph{bi-isomorphism} for short.
Since $\langle\supp_{Q}(\rL)\rangle= Q$ and
$\langle\supp_\Lm(\rL)\rangle= \Lm$,
the maps $\ph_r$ and $\ph_e$ are uniquely determined.
\end{definition}

\begin{remark}
Suppose  that $\ph: \rL \to \rL'$
is  a bi-isomorphism as in Definition \ref{def:bi-isomorphism}.
We did not assume  that $\ph_r$ carries
$\De$ onto $\De'$. However,  this holds automatically,
if $\rL$ and $\rL'$ satisfy (LT5).
In particular, if $\rL$ and $\rL'$ satisfy (LT5) and
$\De' = \De$, then $\ph_r\in \Aut(\De)$.
\end{remark}

\begin{definition}
\label{def:gradingpair}
Suppose that $\rL$ is a
centreless Lie $\Lm$-torus\index{Lie torus!centreless}
of type $\De$.
We set
\[\fg = \rL^0 \andd \fh = \rL_0^0.\]
Then $\fg$ is a finite dimensional simple Lie algebra and $\fh$ is a splitting Cartan subalgebra of $\fg$
\cite{N1,ABFP2}.  Moreover, $\fh$ acts ad-diagonally on $\rL$.   Furthermore,
we can and do identify $\De$ with a root system in the dual space $\fh^*$ in such a way
that $\De(\fg,\fh) = \Dei$,
$\De(\rL,\fh) =  \De$ or $\Dei$,
and, for $\al\in Q$, the root space of
$\rL$ relative to $\al$ is $\rL_\al$ \cite[Prop.~1.2.2]{ABFP2}.
In particular, the type of the split simple Lie algebra $\fg$ is the type
of the root system $\Dei$.
We call $(\fg,\fh)$ the
\emph{grading pair}\index{Lie torus!grading pair}
for
$\rL$.\footnote{In \cite{N1}, the grading pair is defined for
not necessarily centreless Lie tori .  For
centreless Lie tori the two definitions are  equivalent \cite[Prop.~1.2.1]{ABFP2}.}
\end{definition}

\begin{remark}[{\cite[Prop.~3.13]{BN}}]
\label{rem:centroid}
Suppose that
$\rL$ is a $\Lm$-torus, and let $\Ct(\rL)$ be the centroid
of $\rL$.\footnote{Recall that the
\emph{centroid}\index{centroid}
of an algebra
is the associative algebra of all endomorphisms of the algebra that commute
with all left and right multiplications.}
Then, $\Ct(\rL)$ is $\Lm$-graded, with
$\Ct(\rL)^\lm = \set{c\in \Ct(\rL) \suchthat c(\rL^\mu) \subseteq \rL^{\lm+\gm} \text{ for } \mu\in\Lm}$.
Let
\[\Gm(\rL) := \supp_\Lm(\Ct(\rL))\]
\index{centroidal grading group!$\Gm(\rL)$}%
be the support of $\Ct(\rL)$, in which case $\Gm(\rL)$ is a subgroup of $\Lm$
called the
\emph{centroidal grading group}\index{centroidal grading group}
of $\rL$. In that case, $\Ct(\rL)$ is
graded-isomorphic to the  group algebra
$\base[\Gm]$ with its natural
$\Gm$-grading (and hence also with its  $\Lm$-grading).
\end{remark}

\begin{example}[The untwisted centreless Lie torus]
\label{ex:untwisted}
\index{Lie torus!untwisted}
 Let  $\fg$ be a finite dimensional
split simple Lie algebra, let $\fh$ be a splitting Cartan subalgebra of $\fg$,
let $\fg = \oplus_{\al\in \De}  \fg_\al$ be the root space decomposition
of $\fg$ relative to $\fh$, where $\De = \De(\fg,\fh)$.
Set $\rL = \fg \ot \base[\Lm]$ with $Q\times \Lm$-grading given by
$\rL_\al^\lm = \fg_\al \ot \base[\Lm]^\lm$ for $\al\in Q$, $\lm\in \Lm$.
Then $\rL$ is a centreless Lie $\Lm$-torus of type $\De$ which we call the
\emph{untwisted centreless Lie $\Lm$-torus of
type $\De$}.\footnote{The universal central
extension of $\fg \ot \base[\Lm]$
is called the
\emph{toroidal}\index{toroidal Lie algebra}
Lie algebra \cite{MRY}, which is one of the origins
of the term Lie torus.}  The grading pair for $\rL$
is $(\fg\ot 1,\fh\ot 1)$.  The centroid of $\rL$ consists of multiplications
by elements of $\base[\Lm]$,  so $\Gm(\rL) = \Lm$.
\end{example}

In later sections, we will need  the following lemma about bi-isomorphisms
of centerless Lie tori of the same type.

\begin{lemma}
\label{lem:bi-isomorphism}
Let $\rL$ be a Lie $\Lm$-torus of type $\De$ with grading
pair  $(\fg,\fh)$, and let
$\rL'$ be a Lie $\Lm'$-torus of type $\De$ with grading
pair $(\fg',\fh')$.  Let $\Pi =\set{\al_1,\dots,\al_\ell}$
be a base for the root system $\De$, choose
$0\ne e_i\in\fg_{\al_i}$, $0\ne f_i\in\fg_{-\al_i}$ with
$[[e_i,f_i],e_i] = 2e_i$, and
choose $0\ne e'_i\in\fg'_{\al_i}$, $0\ne f'_i\in\fg'_{-\al_i}$ with
$[[e'_i,f'_i],e'_i] = 2e'_i$.
Suppose there is a bi-isomorphism $\ph: \rL \to \rL'$ with $\ph_r$ in the Weyl group
of $\De$.    (Here we are using the notation $\ph_r$ and $\ph_e$ of
Definition \ref{def:bi-isomorphism}.) Then, there exists
a bi-isomorphism $\tilde\ph : \rL \to \rL'$ such that
$\tilde\ph_e = \ph_e$, $\tilde\ph_r = 1$,
$\tilde\ph(e_i) = e_i'$ and $\tilde\ph(f_i) = f_i'$
for all $i$.
\end{lemma}

\begin{proof}  Since $\ph_r$ is in the Weyl
group of $\De$, there exists an automorphism $\eta$ of $\fg$ such
that $\eta(\fh) = \fh$ and $\eta(\fg_\al) = \fg_{\ph_r(\al)}$  for all $\al\in\De$
\cite[Theorem 2(ii), Ch.~VIII, \S~2, no.~2]{B3}.
Moreover, $\eta$ can be
chosen in the form
$\eta = \textstyle\prod_{i}\exp(\ad_\fg(x_i))$, where each
$x_i$ is in a root space of $\fg$ corresponding to a nonzero root.
Thus each $\ad_\rL(x_i)$ is nilpotent and so we may extend
$\eta$ to $\eta = \textstyle\prod_{i}\exp(\ad_\rL(x_i))$.
Then $\eta : \rL \to \rL$ is a bi-isomorphism with $\eta_e = 1$ and $\eta_r = \ph_r$.
So replacing $\ph$ by $\ph\eta^{-1}$, we can assume that $\ph_r = 1$.
Thus, $\ph(\fg_{\al_i}) =  \ph(\fg\cap\rL_{\al_i}) =
\fg'\cap\rL'_{\al_i} = \fg'_{\al_i}$ for all $i$.  So
$\ph(e_i) = a_i e_i'$ and $\ph(f_i) = a_i^{-1} f_i'$ for some
$a_i\in\base^\times$.  Choose a group homomorphism
$\rho : Q \to \base^\times$ such that  $\rho(\al_i) = a_i$ for all $i$,
and define $\tau\in\Aut(\rL)$ by
$\tau(x_\al) = \rho(\al)x_\al$ for $x_\al\in\rL_\al$, $\al\in Q$.
Then $\tau : \rL \to \rL$ is a bi-isomorphism with $\tau_e = 1$ and
$\tau_r = 1$.  So $\ph\tau^{-1}$ is the  required~$\tilde\ph$.
\qed\end{proof}

\section{The construction of EALAs from Lie tori}
\label{sec:construction}
\index{extended affine Lie algebra!construction|(}

In \cite{BGK}, Berman, Gao and Krylyuk gave a  construction
of a family of EALAs starting from a centreless Lie tori of type A$_r$, $r\ge 3$.
In \cite{N2}, Neher simplified this construction and extended it to all types.
Since we will be working with this construction in some detail, we give a careful description
of it in this section. Facts that we note without reference are either straightforward
or  can be found in \cite{N1,N2}.

We first need to establish some notation and assumptions.
Throughout the section, \emph{we assume that
$\rL$ is a centreless Lie $\Lm$-torus of type $\De$  and we let $Q = Q(\De)$}.

Let
\[\Gm = \Gm(\rL) \]
be the centroidal grading  group of $\rL$. We identify
$\Ct(\rL) = \base[\Gm]$ (see Remark \ref{rem:centroid}) by fixing a basis
$\set{t^\gm}_{\gm\in\Gm}$  for $\Ct(\rL)$ satisfying $t^{\gm} t^{\delta} = t^{{\gm}+{\delta}}$.

Next let $\Hom(\Lm,\base)$ be the group
of group homomorphisms from $\Lm$ into $\base$.
Then for  $\theta\in \Hom(\Lm,\base)$, let $\partial_\theta\in \Der(\rL)$ be the \emph{degree derivation}
defined by
\[\partial_\theta(x^\lm) = \theta(\lm)x^\lm\]
for $x^\lm\in \rL^\lm$, $\lm\in \Lm$.  Put
\[\cD = \set{\partial_\theta \suchthat \theta\in \Hom(\Lm,\base)} \andd \CDer(\rL)
= \Ct(\rL) \cD.\]
Then $\CDer(\rL)$ is a subalgebra of the Lie algebra $\Der(\rL)$
with product given by
\begin{equation}
\label{eq:CDermult}
[t^{\gm_1}\partial_{\theta_1}, t^{\gm_2}\partial_{\theta_2}] =
t^{\gm_1+\gm_2}(\theta_1(\gm_2)\partial_{\theta_2}-\theta_2(\gm_1)\partial_{\theta_1}),
\end{equation}
and $\CDer(\rL)$ is
$\Gm$-graded with $\CDer(\rL)^\gm = \Ct(\rL)^\gm\cD$ for $\gm\in \Gm$.

Since $\rL$ is a centreless Lie torus, there is a nondegenerate invariant
$\Lm$-graded form $\form$ on $\rL$, and  that form is unique up to
multiplication by a nonzero scalar \cite[Thm. 2.2 and 7.1]{Y6}.  We fix a choice
of $\form$ on $\rL$.  Then, since $\rL$ is perfect and $\form$ is invariant we have
\[(c(x),y) = (x,c(y))\]
for $x,y\in\rL$, $c\in\Ct(\rL)$.
Now let $\SCDer(\rL)$ be the subalgebra of  $\CDer(\rL)$ consisting
of the derivations in $\CDer(\rL)$ that are skew relative to the form $\form$.
(This subalgebra  does not depend on the choice of the form $\form$.)  Then,
$\SCDer(\rL)$ is a $\Gm$-graded subalgebra of $\CDer(\rL)$ with
\begin{equation}
\label{eq:Cgrad}
\SCDer(\rL)^\gm = t^\gm\set{\partial_\gm\in \cD \suchthat \theta(\gm) = 0}
\end{equation}
for $\gm\in\Gm$. $\SCDer(\rL)$ is called the
\emph{algebra of skew-centroidal derivations of $\rL$}.

Now suppose that $\rD = \oplus_{\gm\in\Gm} \rD_\gm$
is a $\Gm$-graded subalgebra of $\SCDer(\rL)$.
The  \emph{graded-dual space} of $\rD$ is the subspace $\rD^\grd = \oplus_{\gm\in\Gm} (\rD_\gm)^*$ of
$\rD^*$,  where  $(\rD_\gm)^*$ is embedded in
$\rD^*$ by letting its elements act trivially on $\rD_\delta$ for $\delta\ne \gm$
We let
\begin{equation*}
\label{eq:SCDer}
\rC = \rD^\grd,
\end{equation*}
and we give the vector space $\rC$ a $\Gm$-grading by setting
\begin{equation}
\rC_\gm = (\rD_{-\gm})^*
\end{equation}
With this grading, $\rC$ is a $\Gm$-graded $\rD$-module
by means of the contragradient action~$\cdot$ given by
\[(d\cdot f)(e) = - f([d,e]),\]
for $d, e \in \rD$, $f\in \rC$.  We also regard
$\rC$ as  a $\rL$-module with trivial action.
Now define $\sg_\rD : \rL \times \rL \to \rC$ by
\[\sg_\rD(x,y)(d) = (dx|y).\]
Then $\sg_\rD$ is a $\Gm$-graded $2$-cocycle for $L$ with values in
the trivial $\rL$-module  $\rC$.\footnote{In this section and in \S \ref{sec:isotopyEALA}
we use the terminology (cochain, cocycle, coboundary and extension)
from the cohomology of Lie algebras.  See for example
\cite[Chapter I, Exercise 12 for \S 3]{B1}.}
Finally, let $\ev : \Lm \to \rC_0$ denote the
\emph{evaluation map} defined by
\[(\ev(\lm))(\partial_\theta) = \theta(\lm)\]
for $\lm\in\Lm$ and $\partial_\theta \in \rD_0 \subseteq \SCDer(\rL)^0 = \cD$.

With this  background we now present the construction.

\begin{construction}[\cite{N2}] Let $\rL$ be a centreless  Lie $\Lm$-torus of type $\De$.
\label{con:famEALA}
In order to construct an EALA from $\rL$ we need two additional ingredients.  We  assume that
\begin{description}[(4.1a) ]
\item[\ass{a} ]  $\rD$ is a $\Gm$-graded subalgebra of $\SCDer(\rL)$ such that the
map $\ev: \Lm \to \rC_0$ is injective,  where $\Gm=\Gm(\rL)$ and $\rC = \rD^\grd$.
\item[\ass{b}]  $\tau : \rD\times \rD \to \rC$ is a $\Gm$-graded
invariant $2$-cocycle\index{invariant 2-cocycle}
of $\rD$ with values in
the $\rD$-module $\rC$ and $\tau(\rD_0,\rD) = 0$.
\end{description}
Note that the assumption in \ass{b} that $\tau$ is \emph{invariant}  means that
\[\tau(d_1,d_2)(d_3) = \tau(d_2,d_3)(d_1)\]
for $d_i\in \rD$.
Let
\[\rE = \rE(\rL,\rD,\tau)=\rD\oplus \rL\oplus \rC,\]
where $\rC = \rD^\grd$.\footnote{We have changed
the order used in  \cite{N2} of the components
of $\rE$.  This is convenient in the proof of Corollary
\ref{cor:isotopyEALA2} below.}
We  identify $\rD$, $\rL$ and $\rC$ naturally
as subspaces of $\rE$.
Define a product $[\, ,\, ]_E$ on $\rE$ by
\begin{multline*}
[d_1+l_1+f_1,d_2+l_2+f_2]_\rE
= [d_1,d_2] + \left([l_1,l_2]+d_1(l_{2})-d_{2}(l_1)\right)\\
+ \left(d_1\cdot f_{2}-d_{2}\cdot f_1+\sg_{\rD}(l_1,l_{2})+\tau(d_1,d_{2})\right)
\end{multline*}
for $d_i\in \rD$, $l_i\in \rL$, $f_i \in \rC$.
Then,  $\rE$ is a $\Lm$-graded algebra with the direct sum grading.
Next  we extend the bilinear form $\form$ on $\rL$ to a  graded bilinear form
$\form$ on $\rE$ by defining:
\[
(d_1+l_1+f_1\mid d_2+l_2+f_2)=(l_1\mid l_{2})+f_1(d_{2})+f_{2}(d_1).
\]
Finally  let
\begin{equation}
\label{eq:con4}
\rH = \rD_0 \oplus \fh\oplus \rC_0.
\end{equation}
in $\rE$, where $(\fg,\fh) = (\rL^0,\rL_0^0)$ is the grading pair for $\rL$.
Thus we have constructed a triple $((\rE(\rL,\rD,\tau),\form,\rH)$, which we often denote simply
by
$\rE(\rL,\rD,\tau)$\index{extended affine Lie algebra!E(L,D,t)@$\rE(\rL,\rD,\tau)$}.
\end{construction}

Neher has announced the following fundamental result on this construction:

\begin{theorem}[{\cite[Thm. 6]{N2}}]
\label{thm:famEALA} (i) If $\rL$ is a centreless Lie $\Lm$-torus, $\rD$
satisfies \ass{a}  and $\tau$ satisfies \ass{b},
then $\rE(\rL,\rD,\tau)$ is an extended affine Lie algebra
of  nullity $\rank(\Lm)$.

(ii) If $\rE$ is an EALA, then $\rE$  is isomorphic as an EALA to
$\rE(\rL,\rD,\tau)$ for some $\rL$, $D$ and $\tau$ as in (i).
\end{theorem}

\begin{remark}
\label{rem:famEALA}
(i) It  is easy to check that
the EALA $\rE(\rL,\rD,\tau)$ in part (i) of the theorem
does not depend, up to
EALA isomorphism, on the choice of the form $\form$ on $\rL$.

(ii)  Neher actually  states more than we've stated in part (ii)
of the theorem.  Indeed, given an EALA $\rE$,
then $\ccoreE$, with a suitable
grading, is a Lie torus satisfying (LT5) and
$\rE$ is isomorphic as an EALA to $\rE(\ccoreE,\rD,\tau)$
for some $D$ and $\tau$ as in (i) \cite[Thm. 6(ii)]{N2}.

(iii) In particular, we can always choose
$\rL$ in part (ii) of the theorem satisfying (LT5).
In that case, the type of $\rE$ and $\rL$  are the same.

(iv) The Lie torus $\rL$ in part (ii) of the theorem is not uniquely determined
up to bi-isomorphism. However, it is uniquely determined up to
isotopy, a fact that has motivated our work on this topic.
We discuss this in detail in \S \ref{sec:isotopyEALA}  below.
\end{remark}

\begin{remark}
\label{rem:Mext}  Suppose  $\rL$ is a centreless Lie torus, $\Gm = \Gm(\rL)$,
and $\rD$ is a $\Gm$-graded
subalgebra of $\SCDer(\rL)$.  Let $\rM :=\rD\oplus \rL$ and identify $\rD$
and $\rL$ naturally as subspaces of $\rM$.  Then, $\rM$
is a $\Lm$-graded Lie algebra with product $[\, ,\, ]$  given by
\[
[d_1 + l_1,d + l_2]=
[d_1,d_{2}]  + \left([l_1,l_2]+d_1(l_2)-d_2(l_1)\right),
\]
and $\rC$ is an $\rM$-module via $(d +l)\cdot f=d\cdot f$.
The cocycles
$\sg_{\rD}$ and $\tau$ can be extended to cocycles for the Lie algebra
$\rM$ with values in the $\rM$-module $\rC$ by letting
$\sg_{\rD}(\rD,\rM)=\sg_{\rD}(\rM,\rD)=0$ and $\tau(\rL,\rM)=\tau(\rM,\rL)=0$.
If $\rD$ and $\tau$ satisfy \ass{a} and \ass{b}, then
as a  Lie algebra
$\rE(\rL,\rD,\tau)$ is equal to  the
extension $\Ext(\rM,\rC,\sg_\rD+\tau)$ of $\rM$ by $\rC$ using the  cocycle
$\sg_{\rD}+\tau$.
\end{remark}

If $\rL$ is a centreless Lie torus, we let
\[\cP(\rL) = \set{(D,\tau) \suchthat
\text{$D$ satisfies \ass{a} and $\tau$ satisfies \ass{b}}}.\]
We note $\cP(\rL)$ is nonempty
since in particular $(\SCDer(\rL),0)\in \cP(\rL)$.
With this notation Theorem \ref{thm:famEALA}(i) tells us that
Construction \ref{con:famEALA} builds a
family\index{extended affine Lie algebra!family}
\[\set{E(\rL,D,\tau)}_{(D,\tau)\in\cP(\rL)}\]
of EALAs from $\rL$. Also,   Theorem \ref{thm:famEALA}(ii) tells us
that any EALA $\rE$ occurs (up to isomorphism of EALAs) in the family
$\set{E(\rL,D,\tau)}_{(D,\tau)\in\cP(\rL)}$ constructed
from some centreless Lie torus~$\rL$.
\index{extended affine Lie algebra!construction|)}

\section{Isotopy for Lie tori}
\label{sec:isotopy}

In this section, we
recall the definition of isotopy from \cite{ABFP2}.
Throughout the section, we assume that $\rL$ is a Lie $\Lm$-torus of type $\De$
and we let $Q = Q(\De)$.

If $\al\in Q$, we  let
\[\Lm_\al = \Lm_\alpha(\rL) :=  \set{\lm\in\Lm \suchthat \rL^\lm_\al \ne 0}.\]
We call $\Lm_\al$ the \emph{$\Lm$-support of $\rL$ at $\al$}.
Note that $0\in \Lm_\al$ if $0\ne \al\in \Dei$, by (LT2)(ii).

\begin{definition}
\label{def:homotope}
Suppose that $\shom \in \Hom(Q,\Lm)$, where  $\Hom(Q,\Lm)$ is the group
of group homomorphisms from $Q$ into $\Lm$.  We define a new $Q\times \Lm$-graded
Lie algebra
$\rLs$\index{Lie torus!isotope!$\rLs$}
as follows.  As a Lie algebra, $\rLs = \rL$.
The grading on $\rLs$ is given by
\begin{equation}
\label{eq:isograding}
(\rLs)^\lm_\al = \rL_\al^{\lm+\shom(\al)}
\end{equation}
for $\al\in Q$, $\lm\in\Lm$.
\end{definition}

The following necessary and sufficient conditions for $\rLs$ to be a Lie torus are obtained in
\cite[Prop.~2.2.3]{ABFP2}.

\begin{proposition}
\label{prop:isotope}
Let $\shom\in \Hom(Q,\Lm)$, and let $\Pi$ be a base for the root system $\De$.
The following statements are equivalent:
\begin{description}[(a) ]
\item[(a)] $\rLs$ is a Lie torus.
\item[(b)] $\shom(\al) \in \Lm_\al$ for all $0\ne \al\in \Dei$.
\item[(c)] $\shom(\al) \in \Lm_\al$ for all $\al\in \Pi$.
\end{description}
\end{proposition}

\begin{definition}[Isotopes of Lie tori]
\label{def:isotope}
Suppose  that $\shom\in \Hom(Q,\Lm)$.  If $\shom$ satisfies the equivalent conditions in Proposition
\ref{prop:isotope}, we say that $\shom$ is \emph{admissible} for $\rL$.  In that case,
we call the Lie torus
$\rLs$ the
$\shom$-\emph{isotope}\index{Lie torus!isotope}
of $\rL$.
\end{definition}

\begin{remark}
\label{rem:specify}  Suppose that $\Pi = \set{\al_1,\dots,\al_r}$ is a base
for $\De$.  By Proposition~\ref{prop:isotope}, to specify an admissible
$\shom\in \Hom(Q,\Lm)$ for $\rL$, and hence an isotope $\rLs$ of $\rL$,
one can arbitrarily choose $\lm_i\in \Lm_{\al_i}$ for $1\le i \le r$, and then
define $\shom\in \Hom(Q,\Lm)$ with $\shom(\al_i) = \lm_i$ for $1\le i \le r$.
\end{remark}

\begin{definition}
\label{def:isotopy}
Suppose that $\rL$ is a Lie $\Lm$-torus of type $\De$ and
$\rL'$ is a Lie $\Lm'$-torus of type $\De'$.
We say that $\rL$ is
\emph{isotopic}\index{Lie torus!isotopic}
to $\rL'$, written
$\rL\sim \rL'$, if some isotope $\rLs$ of $\rL$ is bi-isomorphic
to $\rL'$.
\end{definition}

Using the facts that $\rL^{(0)} = \rL$, $(\rLs)^{(t)} = \rL^{(\shom + t)}$
and $(\rLs)^{(-\shom)} = \rL$ as $Q\times\Lm$-graded algebras
for $s,t\in \Hom(Q,\Lm)$, it is shown in \cite[\S 2]{ABFP2}
that isotopy is an equivalence relation on the class of Lie tori.

\section{Isotopy in the theory of EALAs}
\label{sec:isotopyEALA}
\index{extended affine Lie algebra!isotopy|(}

As promised at the end of \S \ref{sec:construction}, we now describe the role
that isotopy plays in the structure theory of
EALAs.  Our first theorem on this topic was actually the starting point of our investigation of isotopy.

\begin{theorem}
\label{thm:isotopyEALA1}
Suppose  $\rL$ is a centreless Lie $\Lm$-torus of type $\De$
and $\rL'$ is a centreless Lie $\Lm'$-torus of type $\De'$.
If some member of the family
$\set{E(\rL,D,\tau)}_{(D,\tau)\in\cP(\rL)}$
is isomorphic as an EALA
to some member of the family
$\set{E(\rL',D',\tau')}_{(D',\tau')\in\cP(\rL')}$,
then $\rL$ is isotopic to $\rL'$.
\end{theorem}

\begin{proof}   Let  $(\fg,\fh)$ be the grading pair of
$\rL$.
Since $\De(\rL,\fh) = \De$ or $\Dei$ (see Definition \ref{def:gradingpair}), we can if necessary
replace $\De$ by $\Dei$ so that
\[\De(\rL,\fh) = \De.\]

Suppose now that $\rE = \rE(\rL,\rD,\tau)$, where
$(\rD,\tau)\in\cP(\rL)$.  We use the notation $Q$, $\rC$, $\ev$,
$\rH$ from Construction \ref{con:famEALA}.
Recall that $\rH = \rD_0 \oplus \fh \oplus\rC$, so we can
identify
\[\rH^* = (\rD_0)^* \oplus \fh^* \oplus (\rC_0)^* = \rC_0 \oplus \fh^* \oplus \rD_0.\]
We use $\ev$ to identify $\Lm$ with   a subgroup of $\rC_0 \subseteq \rH^*$.

Let
$R$ be the set of roots of $\rE$ and let $R^\times$ be the
set of nonisotropic roots in $R$.
Then by \cite[\S 5]{N2}, we have
\begin{equation}
\label{eq:T1Rc}
R^\times = \cup_{\al\in\Dec}(\Lm_\al + \al)\subseteq R \subseteq \Lm\oplus Q,
\end{equation}
where  $\Lm_\al = \Lm_\al(\rL)$.  It follows  from this fact and (LT2)(i) that
\begin{equation}
\langle R \rangle = \Lm \oplus Q.
\end{equation}

Now by \cite[\S 5]{N2}, we have
\begin{equation}
\label{eq:T1Eroot}
\rE_{\lm+\al} = \rL_\al^\lm
\end{equation}
for $\lm\in\Lm$, $0\ne \al\in Q$.  It is easy to show using this equation,
\eqref{eq:T1Rc} and (LT3) that the core of $\rE$ (see Definition \ref{def:EALAisom}(ii))
is given by
\begin{equation}
\label{eq:T1Ecore}
\coreE = \rL\oplus \rC.
\end{equation}
(In fact this is implicit in \cite{N2}.)  Hence, since $\rL$ is centreless, we have
\begin{equation}
\label{eq:T1ZEcore}
Z(\coreE) = \rC.
\end{equation}

Next, let $\rE' = \rE(\rL',\rD',\tau')$, where
$(\rD',\tau')\in\cP(\rL')$, and we use all of the above
notation (with primes added) for  $\rL'$ and $\rE'$.

Suppose that $\chi : \rE \to \rE'$ is an
EALA isomorphism.  So
\[\chi(\rH) = \rH',\]
and $\chi$ preserves the forms up multiplication by a nonzero scalar.
It follows that
$\chi(\coreE) = \coreE'$.  So, by \eqref{eq:T1Ecore} and \eqref{eq:T1ZEcore}
(and the corresponding primed equations), we have
\[\chi(\rL\oplus \rC) = \rL'\oplus \rC' \andd \chi(\rC) = \rC'.\]
Thus, there exists a unique Lie algebra isomorphism $\ph : \rL \to \rL'$ with
\begin{equation}
\label{eq:T1chiphi}
\chi(l) \equiv \ph(l) \pmod{C'}
\end{equation}
for $l\in\rL$.  It then follows that
\[\chi^{-1}(l') \equiv \ph^{-1}(l') \pmod{C}\]
for $l'\in\rL'$.

Now, by   \eqref{eq:T1chiphi}, we have
$\ph(\fh) \subseteq (\rH'\oplus \rC')\cap \rL' = \fh'$,
and similarly $\ph^{-1}(\fh') \subseteq \fh$.
So
\begin{equation}
\label{eq:T1phih}
\ph(\fh) = \fh'.
\end{equation}
Hence, if $h\in\fh$,  $\chi(h)-\ph(h) \in \rH'\cap \rC' = \rC_0'$, so
\begin{equation}
\label{eq:T1chiphih}
\chi(h) \equiv \ph(h) \pmod{\rC_0'}.
\end{equation}
Also, $\chi(\rC_0) \subseteq \rH'\cap \rC' = \rC_0'$ and
similarly
$\chi^{-1}(\rC_0') \subseteq \rC_0$.  So
\begin{equation}
\label{eq:T1chiC}
\chi(\rC_0) = \rC_0'.
\end{equation}
But then, by \eqref{eq:T1phih}, \eqref{eq:T1chiphih} and \eqref{eq:T1chiC},
$\chi(\fh\oplus \rC_0) \subseteq \fh'\oplus \rC_0'$; so using the
corresponding property for $\chi^{-1}$, we have
\begin{equation}
\label{eq:T1chifC}
\chi(\fh\oplus \rC_0) = \fh'\oplus \rC_0'.
\end{equation}

Next define $\hat\chi : \rH^* \to {\rH'}^*$ by
\[(\hat\chi(\rho))(\chi(t)) = \rho(t)\]
for $\rho\in\rH^*$, $t\in \rH$.  ($\hat\chi$ is the \emph{inverse dual}
of $\chi|_\rH : \rH \to \rH'$.)  Similarly, define
$\hat\ph : \fh^* \to {\fh'}^*$ by
\[(\hat\ph(\al))(\ph(h)) = \al(h)\]
for $\al\in\fh^*$, $h\in\fh$.  Then, using
$R = \De(\rE,\rH)$ and $\De = \De(\rL,\fh)$ (and the corresponding
primed equations), we see that
\[\hat\chi(R) = R' \andd \hat\ph(\De) = \De'.\]
Thus, since $\langle R \rangle = \Lm \oplus Q$ and $\langle \De \rangle = Q$,
we have
\begin{gather}
\label{eq:T1hat1a}
\hat\chi(\Lm \oplus Q) = \Lm' \oplus Q',\\
\label{eq:T1hat1b}
\hat\ph(Q) = Q'.
\end{gather}

Let $\al\in Q$, and write, using \eqref{eq:T1hat1a},
$\hat\chi(\al) = \lm' + \al'$, where $\lm'\in\Lm'$, $\al'\in Q'$.
Then, $(\lm'+\al')(\chi(h)) = \al(h)$ for $h\in \fh$ (in fact for $h\in\rH$).
Thus, by \eqref{eq:T1chiphih},  $(\lm'+\al')(\ph(h)) = \al(h)$, so
$\al'(\ph(h)) = \al(h)$ for $h\in\fh$.  Hence, $\al'=\hat\ph(\al)$.  Thus, since
$\hat\ph$  is invertible,
\begin{equation}
\label{eq:T1hat2}
\hat\chi(\al) = \shom'(\hat\ph(\al)) + \hat\ph(\al)
\end{equation}
for $\al\in Q$ and some $\shom'\in\Hom(Q',\Lm')$.

Next let $\lm\in \Lm$, and as above write
$\hat\chi(\lm) = \lm' + \al'$, where $\lm'\in\Lm'$, $\al'\in Q'$.
Then, $(\lm'+\al')(\chi(\fh+\rC_0)) = \lm(\fh+\rC_0) = 0$. Hence, by
\eqref{eq:T1chifC}, $(\lm'+\al')(\fh'+\rC_0') = 0$, so $\al'= 0$.
Thus, $\hat\chi(\Lm) \subseteq \Lm'$ and similarly
$\hat\chi^{-1}(\Lm') = \widehat{\chi^{-1}}(\Lm')\subseteq \Lm$.  So
\begin{equation}
\label{eq:T1hat3}
\hat\chi(\Lm) = \Lm'.
\end{equation}

Using $\eqref{eq:T1hat1b}$ and $\eqref{eq:T1hat3}$, we let
$\ph_r = \hat\ph|_Q : Q \to Q'$ and $\ph_e = \hat\chi|_\Lm : \Lm \to \Lm'$.
Then, using \eqref{eq:T1Eroot} and \eqref{eq:T1hat2}, we have, for $0\ne \al\in Q$, $\lm\in\Lm$,
\begin{align*}
\chi(\rL_\al^\lm) &= \chi(\rE_{\lm+\al})
= \rE'_{\hat\chi(\lm+\al)}
= \rE'_{\hat\chi(\lm)+\shom'(\hat\ph(\al)) + \hat\ph(\al)}\\
&= (\rL')^{\hat\chi(\lm)+\shom'(\hat\ph(\al))}_{\hat\ph(\al)}
= ((\rL')^\sip)^{\ph_e(\lm)}_{\ph_r(\al)}.
\end{align*}
So by \eqref{eq:T1chiphi}, we have
\begin{equation*}
\label{eq:T1equiv}
\ph(\rL_\al^\lm) = ((\rL')^\sip)^{\ph_e(\lm)}_{\ph_r(\al)}
\end{equation*}
for $0\ne \al\in Q$, $\lm\in\Lm$.  But then, by (LT3), this equation holds for \emph{all}
$\al\in Q$, $\lm\in\Lm$.  Therefore, since $\rL$ is a Lie $\Lm$-torus of type $\De$,
so is $(\rL')^\sip$.  Hence  $\shom'$ is admissible for $\rL'$, and
$\ph$ is a bi-isomorphism of $\rL$ onto $(\rL')^\sip$.
\qed\end{proof}

We next consider the relationship between the family of EALAs constructed from a
Lie torus $\rL$ and the family of EALAs constructed from an isotope of $\rL$.
For this purpose, we assume now  that
$\rL$ is a centreless Lie $\Lm$-torus of type $\De$,
$Q = Q(\De)$, $\shom\in\Hom(Q,\Lm)$ is admissible
for $\rL$, and $\rLs$ is the $\shom$-isotope of~$\rL$.

We fix a nondegenerate $\Lm$-graded invariant symmetric bilinear
form $\form$ on $\rL$ and we use that same form on $\rLs$.

Let $\partial_{\theta}^\si$  be the degree derivation
of $\rLs$ determined by $\theta\in \Hom(\Lambda,\base)$.
Define  the $\Lm$-graded $\Ct(\rL)$-module isomorphism
\[\Psi:\CDer(\rL)\rightarrow \CDer(\rLs) \quad \text{by}
\quad \Psi(c\partial_{\theta})=c\partial_{\theta}^\si.\]
Using \eqref{eq:CDermult} we see that $\Psi$ is also a
Lie algebra isomorphism.
Let $h_{\theta}\in\mathfrak{h}$ be given by
\[\alpha(h_{\theta})=\theta(s(\alpha))\]
for $\alpha\in Q$, and define the $\Lambda$-graded  $\Ct(\rL$)-module
monomorphism
\[\Omega:\CDer(\rL)\rightarrow \rL \quad\text{by}\quad\Omega(c\partial_{\theta})=ch_{\theta}\]
for  $c\in \Ct(\rL)$.
For $l\in(\rLs)_{\alpha}^{\lambda}=\rL_{\alpha}^{\lambda+s(\alpha)}$,
we have $\partial_{\theta}(l)=\theta(\lambda+s(\alpha))l=\partial_{\theta}^\si(l)+[h_{\theta},l]$.
Thus,
\[\partial_{\theta}^\si=\partial_{\theta}-\ad_\rL(h_{\theta}) \andd \Psi(d)=d-\ad_\rL(\Omega(d))\]
for $d\in \CDer(\rL)$.
Since the derivations in  $\ad_\rL(\rL)$ are skew, we see that
$\Psi$ maps $\SCDer(\rL)$ into $\SCDer(\rLs)$.  We also note  using \eqref{eq:CDermult} that
\begin{align*}
\Omega([t^{\gamma_1}\partial_{\theta_1},t^{\gamma_{2}}\partial_{\theta
_{2}}])  & =t^{\gamma_1}t^{\gamma_{2}}(\theta_1(\gamma_{2})h_{\theta_{2}%
}-\theta_{2}(\gamma_1)h_{\theta_1})\\ & =t^{\gamma_1}\partial_{\theta_1}(t^{\gamma_{2}}h_{\theta_{2}}%
)-t^{\gamma_{2}}\partial_{\theta_{2}}(t^{\gamma_1}h_{\theta_1})\\
& =t^{\gamma_1}\partial_{\theta_1}(\Omega(t^{\gamma_{2}}\partial
_{\theta_{2}}))-t^{\gamma_{2}}\partial_{\theta_{2}}(\Omega(t^{\gamma_1%
}\partial_{\theta_1})),
\end{align*}
for $\gm_i\in\Gm$ and $\theta_i\in\Hom(\Lm,\base)$,
so $\Omega$ is a derivation of $\CDer(\rL)$ into  $\rL$.

If $(\rD,\tau)\in\cP(\rL)$, let
\begin{equation}
\label{eq:Ds}
\rD^\si=\Psi(\rD),
\end{equation}
in which case
$\rD^\si$ is a $\Gm$-graded  subalgebra of  $\SCDer(\rL)$.  Also let
\[\psi:\rD\rightarrow \rD^\si \andd \omega:\rD\rightarrow \rL\]
be the restrictions of
$\Psi$ and $\Omega$ to $\rD$. Then $\psi$ is a $\Lm$-graded Lie algebra
isomorphism, $\omega$ is a $\Lm$-graded vector space monomorphism, and
\begin{equation}
\label{eq:psiomega}
\psi(d)=d-\ad_\rL(\omega(d))
\end{equation}
for $d\in\rD$.   Finally let
\[\rC^\si=(\rD^\si)^\grd,\]
which  is $\Gm$-graded as in \eqref{eq:Cgrad}, and define a $\Lm$-graded
cochain
\begin{equation}
\label{eq:taus}
\tau^\si:\rD^\si\times
\rD^\si\rightarrow \rC^\si \quad\text{by}\quad
\tau^\si(\psi d_1,\psi d_{2})(\psi d_{3})=\tau(d_1,d_{2})(d_{3})
\end{equation}
for $d_{i}\in \rD$.

\begin{theorem}
\label{thm:isotopyEALA2}
Suppose that $\rL$ is a centreless Lie $\Lm$-torus, $s\in \Hom(Q,\Lambda)$ is admissible for $\rL$,
$(\rD,\tau)\in\cP(\rL)$, and $(\rD^\si,\tau^\si)$ is defined by \eqref{eq:Ds} and \eqref{eq:taus}.
Then $(\rD^\si,\tau^\si) \in \cP(\rLs)$.
Moreover, $\rE = \rE(\rL,\rD,\tau)$ and $\rE^\si =\rE(\rLs,\rD^\si,\tau^\si)$
are isomorphic as EALAs. To be more precise,
define $\chi : \rE \to\rE^\si$ by
\[\chi(x) =
\psi(d) + \left(\omega(d)+l\right) + \eta(x)
\]
for $x = d + l + c\in\rE$, $d\in\rD$, $l\in\rL$, $f\in\rC$, where $\eta:\rE\rightarrow \rC^\si$ is given by
\[
\eta(x)(\psi(d'))=f(d')-(l+\frac{1}{2}\omega(d)\mid
\omega(d'))
\]
for $d'\in\rD$.  Then $\chi$
is an isometry and a  $\Lm$-graded
Lie algebra isomorphism that maps
$H$ to $H^\si=\rD_{0}^\si\oplus\mathfrak{h}\oplus \rC_{0}^\si$.
\end{theorem}

\begin{proof}
Define  $\hat{\psi}:\rC\rightarrow \rC^\si$ by $\hat
{\psi}(f)(\psi(d))=f(d)$ for $d\in D$, $f\in \rC$.  Clearly, $\hat{\psi}$ is a
$\Gm$-graded  vector space isomorphism.
Let $\ev^\si$ be the evaluation map for $\rLs$.
Then, for $\lambda\in\Lambda$ and $\partial_{\theta}^\si\in D_{0}^\si$, we have
\[
(\hat{\psi}(\ev(\lambda)))(\partial_{\theta}^\si)=(\hat{\psi}(\ev(\lambda
)))(\psi(\partial_{\theta}))=(\ev(\lambda))(\partial_{\theta})=\theta(\lambda),
\]
so  $\hat{\psi}\mid_{\rC_0}\circ \ev=\ev^\si$.  Hence, $\ev^\si$ is injective; i.e.,
$D^\si$ satisfies (4.1a).

It will be convenient to view $\rE$ as $\Ext(\rM,\rC,\sg_{D}+\tau)$, the extension
of $\rM=D\oplus \rL$ by $\rC$ using the cocycle $\sg_{D}+\tau$ as in Remark
\ref{rem:Mext}.  Let $\rM^\si=D^\si\oplus \rL$ and set
\[\xi(d+l)=\psi(d)+\omega(d)+l\]
for $d\in D$, $l\in \rL$.  We claim
\begin{equation}
\xi:\rM\rightarrow \rM^\si\text{ is a Lie algebra isomorphism.}\label{eq:xiisom}
\end{equation}
Clearly, $\xi$ is bijective since $\psi$ is bijective.  Since $\psi
:D\rightarrow D^\si$ is an isomorphism and $\rL$ is an ideal in $\rM^\si$, the
map $d+l\rightarrow\psi(d)+\rL=\xi(d+l)+\rL$ is a homomorphism $\rM\rightarrow
\rM^\si/\rL$.  Thus, for $m_1,m_2\in \rM$,
\[
z:=\xi([m_1,m_2])-[\xi(m_1),\xi(m_2)]\in \rL.
\]
On the other hand, $\zeta:\rM\rightarrow \Der(\rL)$ given by
$\zeta(m)=\ad_{\rM}(m)\mid_{\rL}$ and similarly $\zeta^\si:\rM^\si\rightarrow \Der(\rL)$ are
homomorphisms.  Moreover,
\[
\zeta^\si(\xi(d+l))=\psi(d)+\ad_\rL(\omega(d))+\ad_\rL(l)=d+\ad_\rL(l)=\zeta(d+l),
\]
so $\zeta^\si\circ\xi=\zeta$.  Thus,
\[
\ad_\rL(z)=\zeta^\si(z)=\zeta([m_1,m_2])-[\zeta(m_1),\zeta (m_2)]=0.
\]
Since $\rL$ is centreless, $z=0$ and $\xi$ is an isomorphism.

We also note that for $d,d'\in D$, $l\in \rL$, $f\in \rC$,
\begin{align*}
(\xi(d+l)\cdot\hat{\psi}(f))(\psi(d'))  & =(\psi(d)\cdot\hat{\psi
}(f))(\psi(d'))
=-\hat{\psi}(f)([\psi(d),\psi(d')])\\
& =-f([d,d'])=(d\cdot f)(d')
=\hat{\psi}((d+l)\cdot f)(\psi(d'),
\end{align*}
so
\begin{equation}
\xi(m)\cdot\hat{\psi}(f)=\hat{\psi}(m\cdot f)\label{eq:actionisom}
\end{equation}
for $m\in \rM$, $f\in \rC$.

Recall that the group $Z^{2}(\rM,\rC)$ of $2$-cocycles on $\rM$  with values
in $\rC$ consists of all
alternating bilinear maps $\mu:\rM\times \rM\rightarrow \rC$  such that
\[
\sum\nolimits_{(i,j,k)\circlearrowleft}
\left(\mu([m_{i},m_{j}],m_{k})-m_{i}\cdot\mu(m_{j},m_{k})\right)=0
\]
for  $m_1,m_2,m_{3}\in \rM$ where $(i,j,k)\circlearrowleft$ means $(i,j,k)$
is a cyclic permutation of $(1,2,3)$.  We can use $\xi$ and $\hat{\psi}$ to
transfer $\mu$ to $\tilde{\mu}:\rM^\si\times \rM^\si\rightarrow \rC^\si$ with
\[\tilde{\mu}(\xi(m_1),\xi(m_2))=\hat{\psi}(\mu(m_1,m_2)).\]
Using
(\ref{eq:xiisom}) and (\ref{eq:actionisom}), we see that $\rho(\mu)=\tilde
{\mu}$ defines a group isomorphism
$\rho:Z^{2}(\rM,\rC)\rightarrow Z^{2} (\rM^\si,\rC^\si)$.   Furthermore,
the map
$\xi\oplus\hat{\psi}:\Ext(\rM,\rC,\mu)\rightarrow
\Ext(\rM^\si,\rC^\si,\rho(\mu))$  is a Lie algebra isomorphism.

Let $\tilde{\tau}=\rho(\tau)$.   Then
$\tilde{\tau}(\rM^\si,\rL)=\tilde{\tau}(\xi(\rM),\xi(\rL))=\hat{\psi}(\tau(\rM,\rL))=0$
and thus also $\tilde{\tau}(\rL,\rM^\si)=0$.  Moreover,  for $d_1,d_2,d_{3}\in
D$, we have
\[
\tilde{\tau}(\psi(d_1),\psi(d_2))(\psi(d_{3}))=\hat{\psi}(\tau(d_1
,d_2))(\psi(d_{3}))=\tau(d_1,d_2)(d_{3}),
\]
so $\tilde{\tau}$ is the extension of $\tau^\si$ to $\rM^\si$ and
$\tau^\si$ is a cocycle.  The invariance of $\tau^\si$ is immediate from
the invariance of $\tau$, while the equation $\tau^\si(D_{0}^\si
,D^\si)=0$ is immediate from the definition of $\tau^\si$ and $\psi
(D_{0})=D_{0}^\si$.  This shows
$(D^\si,\tau^\si)\in\mathcal{P}
(\rL^\si)$.

It remains to prove the statement about $\chi$.
To do this it will be helpful for us to write an element
$x = d + l + f \in \rE$, where $d\in \rD$, $l\in \rL$, $f\in C$,
as a column vector
\[x=\left[ \begin{array}[c]{c} d\\ l\\f \end{array}\right].\]
Similarly we express elements of $\rE^\si$ as column vectors.
This notation allows us to write linear maps (for example from
$\rE$ to $\rE^\si$) as matrices.

Define
\[\omega^{\#}:\rL\rightarrow \rC \quad\text{by} \quad \omega
^{\#}(l)(d)=(l\mid\omega(d)).\]
We see that for $x = d + l + c\in\rE$, $d\in\rD$, $l\in\rL$, $f\in\rC$
and for $d'\in \rD$, we have
\[\eta(x)(\psi(d'))   =(f-\omega^{\#}(l+\frac{1}{2}\omega
(d)))(d')\\
=\hat{\psi}(f-\omega^{\#}(l+\frac{1}{2}\omega(d)))(\psi(d')),
\]
so $\eta(x)=\hat{\psi}(f)-(\hat{\psi}\circ\omega^{\#})(l)-\frac{1}{2}
(\hat{\psi}\circ\omega^{\#}\circ\omega)(d)$.
Hence, in matrix form
\[
\chi=\left[
\begin{array}
[c]{ccc}%
\psi & 0 & 0\\
\omega & \Id & 0\\
-\frac{1}{2}\hat{\psi}\circ\omega^{\#}\circ\omega & -\hat{\psi}\circ
\omega^{\#} & \hat{\psi}%
\end{array}
\right]  .
\]

Let
\[
\nu=\left[
\begin{array}
[c]{ccc}%
0 & 0 & 0\\
\omega & 0 & 0\\
0 & -\omega^{\#} & 0
\end{array}
\right]  ,
\]
so $\nu:\rE\rightarrow \rE$.  The definition of $\omega^{\#}$  shows that $\nu$
is a skew transformation of $\rE$.  Thus, $\exp(\nu)$ is an isometry of $\rE$.
We can write
\[
\chi=\left[
\begin{array}
[c]{ccc}%
\psi & 0 & 0\\
0 & \Id & 0\\
0 & 0 & \hat{\psi}%
\end{array}
\right]  \left[
\begin{array}
[c]{ccc}%
\Id & 0 & 0\\
\omega & \Id & 0\\
-\frac{1}{2}\omega^{\#}\circ\omega & -\omega^{\#} & \Id
\end{array}
\right]  =\left[
\begin{array}
[c]{ccc}%
\psi & 0 & 0\\
0 & \Id & 0\\
0 & 0 & \hat{\psi}%
\end{array}
\right]  \exp(\nu).
\]
The first factor is an isometry $\rE\rightarrow \rE^\si$ by the definition of
$\hat{\psi}$.  Thus, $\chi:\rE\rightarrow \rE^\si$ is an isometry.

Since $\psi$ and $\omega$ are $\Lambda$-graded, so are $\hat\psi$,
$\omega^\#$, $\eta$ and $\chi$.
Also, $\omega
(\rD_{0})\subseteq\mathfrak{h}$, so $\chi(H)\subseteq\psi(\rD_{0})\oplus
(\omega(\rD_{0})+\mathfrak{h)}\oplus\eta(\rE_{0})\subseteq H^\si$.
But since $\rD_0^\si = \psi(\rD_0)$ and
 $\rC_0^\si = \hat\psi(\rC_0)$, we  have $\dim(\rH^\si) = \dim (\rH)$,
so $\chi(\rH) = \rH^\si$.

It remains  to show that $\chi$ is a homomorphism of Lie algebras.  Let
$\sg'=\rho^{-1}(\sg_{D^\si})\in Z^{2}(\rM,\rC)$.
Then,  $\sg_{D^\si}(\xi(m_1),\xi(m_2)) = \hat\psi(\sg'(m_1,m_2))$ for
$m_1,m_2\in\rM$, so
\[\sg'(m_1,m_2)(d) = \sg_{D^\si}(\xi(m_1),\xi(m_2))(\psi(d))\]
for $d\in\rD$.
We  will   show that $\sg'=\sg_{\rD}+\delta(\kappa)$ where
$\kappa:\rM\rightarrow \rC$ is the $1$-cocycle%
\[
\kappa(d+l)=\frac{1}{2}\omega^{\#}(\omega(d))+\omega^{\#}(l)
\]
and $\delta(\kappa)$ is the coboundary of  $\kappa$; i.e.,
\[
\delta(\kappa)(m_1,m_{2})=m_{2}\cdot\kappa(m_1)-m_1\cdot\kappa(m_{2}%
)-\kappa([m_1,m_{2}])\text{.}
\]
Indeed, using linearity and skew symmetry, it suffices to check this for the pairs
$(l_1,l_{2})$, $(l_1,d_{2})$, and $(d_1,d_{2})$.  For this recall that
$\rL\cdot \rC=0$, $[\omega(\rD),\omega(\rD)]\subseteq\lbrack \Ct(\rL)\mathfrak{h,}\Ct(\rL)\mathfrak{h}%
]=0$, and that $\omega$ is a derivation of $\rD$ into $\rL$.  We have
\begin{align*}
\sg'(l_1,l_{2})(d)  & =\sg_{\rD^\si}(l_1,l_{2})(\psi(d))
 =(\psi(d)(l_1)\mid l_{2})\\
& =((d-\ad_\rL(\omega(d))(l_1)\mid l_{2})
=(d(l_1)-[\omega(d),l_1]\mid l_{2})\\
& =(d(l_1)\mid l_{2})-(\omega(d)\mid\lbrack l_1,l_{2}])
=(\sg_{\rD}(l_1,l_{2})-\omega^{\#}([l_1,l_{2}]))(d),
\end{align*}
while
\[\delta(\kappa)(l_1,l_{2})
=l_1\cdot\kappa(l_{2})-l_{2}\cdot\kappa(l_1)-\kappa([l_1,l_{2}])
=-\kappa([l_1,l_{2}])=-\omega^{\#}([l_1,l_{2}])).\]
Also,
\begin{align*}
\sg'(d_1,d_{2})(d)  & =\sg_{\rD^\si}(\psi(d_1)+\omega
(d_1),\psi(d_{2})+\omega(d_{2}))(\psi(d))\\
& =(\psi(d)(\omega(d_1))\mid\omega(d_{2}))\\
& =((d-\ad_\rL(\omega(d))(\omega(d_1))\mid\omega(d_{2}))\\
& =(d(\omega(d_1))\mid\omega(d_{2}))
\end{align*}
and
\begin{align*}
2(\delta(\kappa)&(d_1,d_{2}))(d)
=(d_1\cdot\omega^{\#}(\omega
    (d_{2}))-d_{2}\cdot\omega^{\#}(\omega(d_1))-\omega^{\#}(\omega([d_1,d_{2}])))(d)\\
& =-(\omega(d_{2})\mid\omega([d_1,d]))+(\omega(d_1)\mid\omega([d_{2},d]))
     \ -(\omega([d_1,d_{2}]))\mid\omega(d))\\
& =-(\omega(d_{2})\mid d_1(\omega(d))-d(\omega(d_1)))
    +(\omega(d_1)\mid d_{2}(\omega(d))-d(\omega(d_{2})))\\
& \qquad\qquad\qquad\qquad -(d_1(\omega(d_{2}))-d_{2}(\omega(d_1))\mid\omega(d))\\
& =(\omega(d_{2})\mid d(\omega(d_1)))-(\omega(d_1)\mid d(\omega(d_{2})))\\
& =2(d(\omega(d_1))\mid\omega(d_{2})).
\end{align*}
Finally,
\[\sg'(d_1,l_{2})(d)   =\sg_{\rD^\si}(\psi(d_1)+\omega
(d_1),l_{2})(d) =(d\cdot\omega(d_1),l_{2})\]
and
\begin{align*}
(\delta(\kappa)(d_1,l_{2}))(d)
& =(d_1\cdot\omega^{\#}(l_{2})-l_{2}\cdot\frac{1}{2}\omega^{\#}(\omega(d_1))
-\omega^{\#}([d_1,l_{2}]))(d)\\
& =-(l_{2}\mid\omega([d_1,d]))+(l_{2}\mid d_1\cdot\omega(d))\\
& =(l_{2}\mid d\cdot\omega(d_1)).
\end{align*}
So  $\sg'=\sg_{\rD}+\delta(\kappa)$ as claimed.

Since  $\sg'+\tau=\sg_{D}+\tau+\delta(\kappa)$, the map
$m+f\rightarrow m+f-\kappa(m)$ for $m\in \rM$, $f\in \rC$ is a Lie algebra
isomorphism of $\rE=\rE(\rL,D,\tau)$ with $\Ext(\rM,\rC,\sg'+\tau)$.  We can
write this isomorphism as
\[
\left[
\begin{array}
[c]{ccc}
\Id & 0 & 0\\
0 & \Id & 0\\
-\frac{1}{2}\omega^{\#}\circ\omega & -\omega^{\#} & \Id
\end{array}
\right]  :\left[
\begin{array}
[c]{c}
D\\
\rL\\
\rC
\end{array}
\right]  \rightarrow\left[
\begin{array}
[c]{c}
D\\
\rL\\
\rC
\end{array}
\right]
\]
Since $\sg_{D^\si}+\tau^\si=\rho(\sg'+\tau)$, we also have
that
\[
\xi\oplus\hat{\psi}=\left[
\begin{array}
[c]{ccc}
\psi & 0 & 0\\
\omega & \Id & 0\\
0 & 0 & \hat{\psi}
\end{array}
\right]
\]
is a Lie algebra isomorphism of $\Ext(\rM,\rC,\sg'+\tau)$ with
$\rE^\si$.  Thus,
\[
\chi=\left[
\begin{array}
[c]{ccc}
\psi & 0 & 0\\
\omega & \Id & 0\\
0 & 0 & \hat{\psi}
\end{array}
\right]  \left[
\begin{array}
[c]{ccc}
\Id & 0 & 0\\
0 & \Id & 0\\
-\frac{1}{2}\omega^{\#}\circ\omega & -\omega^{\#} & \Id
\end{array}
\right]
\]
is a Lie algebra isomorphism from $\rE$ to  $\rE^\si$.
\qed\end{proof}

\begin{corollary}
\label{cor:isotopyEALA2}
Let $\rL$ be a centreless Lie $\Lm$-torus of type $\De$,
and let $\rL'$ be a centreless Lie $\Lm'$-torus of type $\De'$.
Suppose that $\rL$ is isotopic to $\rL'$. Then, there is a bijection
$(D,\tau) \mapsto (D',\tau')$ from $\cP(\rL)$ onto $\cP(\rL')$ such that
$E(\rL,D,\tau)$ is isomorphic as an EALA to $E(\rL',D',\tau')$
for all $(D,\tau)\in \cP(\rL)$.
\end{corollary}

\begin{proof}  If $\rL$ and $\rL$ are bi-isomorphic,  the result is clear.
So we can assume that $\rL' = \rLs$ is an isotope of $\rL$.  Define
$T : \cP(\rL) \to \cP(\rL')$ by $T(\cD,\tau) = (\cD^\si,\tau^\si)$
using the notation of  Theorem \ref{thm:isotopyEALA2}.  Similarly, we have a map
$T' :  \cP(\rL') \to \cP({\rL'}^{(-\shom)}) = \cP(\rL)$, and it is
clear that $T$ and $T'$ are inverses of one another.  So $T$ is a bijection.
\qed\end{proof}

Less precisely (but more succinctly), Corollary \ref{cor:isotopyEALA2}
says that  if the centerless Lie tori $\rL$ and $\rL'$ are isotopic then the families
$\set{E(\rL,D,\tau)}_{(D,\tau)\in\cP(\rL)}$ and
$\set{E(\rL',D',\tau')}_{(D',\tau')\in\cP(\rL')}$ of EALAs are
\emph{bijectively isomorphic}\index{extended affine Lie algebra!family!bijectively isomorphic}.
Using this language, Theorem \ref{thm:famEALA},
Theorem \ref{thm:isotopyEALA1} and Corollary \ref{cor:isotopyEALA2}
together tell us that Construction \ref{con:famEALA} provides a
one-to-one correspondence between centreless Lie tori up to isotopy
and families of EALAs up to bijective  isomorphism.

\begin{remark}
\label{rem:EALAexplicit}  Suppose that $\rL$ is a centreless Lie torus.
To construct all of the members of the family $\set{E(\rL,D,\tau)}_{(D,\tau)\in\cP(\rL)}$
of EALAs  corresponding to $\rL$, one must do the following:
(i)  Calculate the centroid $\Ct(\rL)$ and identify  it explicitly with $\base[\Gm]$ (in which case
$\SCDer(\rL)$ is completely understood from  \eqref{eq:Cgrad});
(ii) Determine the $\Gm$-graded subalgebras $D$ of $\SCDer(\rL)$ which satisfy
\ass{a};
(iii) For each $\rD$ in (ii), determine all 2-cocycles of $D$
satisfying \ass{b}.  We note that (i) can be carried out for each type $\De$
using the coordinatization theorems (see \S \ref{sec:coordinates} below)
and the results and methods in \cite[\S 5]{BN}.
However, finding   a general approach to (ii) and (iii)
seems to be much more difficult. (See \cite[p.3]{BiN} for some remarks on problem (iii) and
see \cite[Remark 3.71(b)]{BGK} for an example of a nontrivial choice of $\tau$.)
We do however note that the  tasks (ii) and (iii) are independent of the type $\De$; in fact they depend
only on the rank of  $\Gm$.
\end{remark}
\index{extended affine Lie algebra!isotopy|)}

\section{Coordinatization of Lie tori}
\label{sec:coordinates}
\index{Lie torus!coordinatization theorems|(}

In the classical study of finite dimensional split simple Lie  algebras,
one can construct the algebras using special linear, orthogonal and symplectic
matrix constructions along with some exceptional constructions
starting from Jordan algebras and alternative algebras.
(See for example \cite[\S III.8, IV.3 and IV.4]{S}.)

Similar  \emph{coordinatization theorems} have been proved
for centreless Lie tori. Very roughly speaking these results show that
a centerless Lie torus of a given type  can be constructed
as a ``matrix algebra'' over a ``coordinate torus'', which
for some types is nonassociative.
References for each type are:
A$_1$: \cite{Y1}; A$_2$:   \cite{BGKN,Y4};
A$_l (l \ge 3)$, $D_l (l\ge 4)$, $E_6$, $E_7$ and $E_8$: \cite{BGK,Y2};
$\text{B}_2 = \text{C}_2$: \cite{AG,BY}; B$_l (l\ge 3)$: \cite{AG,Y5,AB};
C$_l (l\ge 3)$: \cite{AG,BY};
F$_4$ and G$_2$: \cite{AG};
BC$_1$ \cite{AFY}; BC$_2$: \cite{F1}; and finally BC$_l (l\ge 3)$: \cite{AB}.
We note that these theorems were in some cases proved before the notion
of Lie torus had been introduced and were instead formulated
in the language of EALAs.  Also, in some cases the theorems were
proved over the complex field.  Nevertheless the translation
to the language of Lie tori and the extension  to
arbitrary base fields of characteristic 0 are not difficult.

In the following sections, we recall the coordinatization theorems for types
A$_1$, A$_2$, A$_l (l\ge 3)$ and C$_l (l\ge 4)$,   which use
respectively Jordan tori, alternative tori, associative
tori and associative tori with involution as coordinates.
First we recall  the definitions of these types of tori.

Recall that a \emph{Jordan algebra} is an algebra with commutative product satisfying
$a^2(ab) = a(a^2 b)$, and an \emph{alternative algebra} is an algebra with product satisfying
$a(ab) = a^2b$ and $(ba)a = ba^2$ (see for example \cite{S} or \cite{Mc2}).

\begin{definition}
\label{def:tori}  Suppose that $\rA =  \bigoplus_{\lm\in\Lm}\rA^\lm$
is a $\Lm$-graded Jordan algebra, alternative algebra
or associative algebra.  $\rA$ is said to be a
\emph{Jordan, alternative or associative $\Lm$-torus}
\index{Jordan torus}
\index{alternative torus}
\index{associative torus}
respectively if:
\begin{description}[(ii) ]
\item[(i)] For $\lm \in \supp_\Lm(\rA)$, $\rA^\lm$ is spanned by a single
invertible  element
of $\rA$.\footnote{We are using the standard definitions
of invertibility for Jordan algebras \cite[p.210]{Mc2},
alternative algebras \cite[p.38]{S} and associative algebras.}
\item[(ii)] $\Lm$ is generated as a group by $\supp_\Lm(\rA)$.
\end{description}
Note of course that any associative torus is an alternative torus.
An
\emph{associative $\Lm$-torus with involution}\index{associative torus with involution}
is a pair $(\rA,\iota)$, where
$\rA$ is an associative $\Lm$-torus and $\iota$ is a graded involution
of $\rA$.
\end{definition}

\begin{remark}
\label{rem:supporttori}  Suppose that $\rA$ is a Jordan or alternative  $\Lm$-torus,
and let $S = \supp_\Lm(\rA)$. Then, it is clear that  $-S= S$.
Moveover, if $\rA$ is alternative, it is easy to see that
$S$ is closed under addition, so $S = \Lm$.  However,
this is not true in the Jordan case, although $S$
is closed under the operation $(\lm,\mu)\mapsto \lm + 2\mu$.
\end{remark}

\begin{example}
\label{ex:quantum}  Let $\bq = (q_{ij})$ be an $n\times n$-matrix
over $\base$ with $q_{ii} = 1$ and $q_{ij} = q^{-1}_{ji}$ for all $i,j$.
Let
$\base_\bq = \base_\bq[x_1^{\pm1},\dots,x_n^{\pm1}]$ be  the associative algebra
over $\base$  presented by the generators $x_i, x^{-1}_i$,  $i=1, \ldots, n$,
subject to the relations
$$x_ix_i^{-1} = x_i^{-1}x_i =1 \andd x_j x_i = q_{ij} x_i x_j, \quad 1\leq i, j \leq n.$$
Then, $\base_\bq$ has a natural $\Zn$-grading with
$(\base_\bq)^{(l_1,\dots,l_n)} = \base x_1^{l_1}\dots x_n^{l_n}$.
The $\Zn$-graded algebra $\base_\bq$ is an
associative  $\Zn$-torus\index{associative torus}
called the
\emph{quantum torus}\index{quantum torus}
determined by $\bq$.\footnote{If $n=0$, we interpret $\bq=\emptyset$ and $\base_\bq = \base$.}
Conversely, it is easy to see that any associative torus
is isograded-isomorphic to a quantum torus.  Recently,  \emph{rational quantum tori},
that is quantum tori $\base_\bq$ with each $q_{ij}$ a root of unity, have
been classified up to isograded-isomorphism by \text{K.-H.}~Neeb in \cite[Thm.~ III.4]{Ne}
as tensor products of rational quantum $\mathbb Z^2$-tori and an algebra of Laurent polynomials.  (In \cite[Thm.~4.8]{H}, J.T.~Hartwig
gave a tensor product decomposition of rational quantum tori but he did not give  a condition for
isomorphism.) \end{example}

In the following sections, we will give some further examples of the tori defined in Definition
\ref{def:tori}.
\index{Lie torus!coordinatization theorems|)}

\section{Type A$_1$}
\label{sec:A1}

In this section we assume that
\[\De = \set{\al,0,-\al}\]
is the root system of type A$_1$, and $Q = Q(\De) = \bbZ \al$.

To construct a centreless Lie $\Lm$-torus of type $\De$,
we use the Tits-Kantor-Koecher  construction.
(See \cite[p.13]{Mc2} for the facts mentioned below about this construction.)
We begin with a Jordan $\Lm$-torus $\rA$ and
we define operators
$V_{x,y}\in\End(\rA)$ for
$x,y\in\rA$ by $V_{x,y} z = \{x,y,z\} := 2((xy)z + (zy)x-(zx)y)$.
These operators satisfy the identity $[V_{x,y},V_{z,w}] = V_{\{x,y,z\},w} - V_{z,\{y,x,w\}}$,
so their span $V_{\rA,\rA}$ is a Lie algebra under the commutator
product.  Moreover, the map $*$ defined by $V_{x,y} \mapsto -V_{y,x}$ is a well-defined
automorphism of $V_{\rA,\rA}$.
We now set
\[\rL(\rA) = \TKK(\rA) := \rA_1 \oplus V_{\rA,\rA}\oplus  \rA_{-1},\]
\index{Tits-Kantor-Koecher Lie algebra!$\TKK(\rA)$}%
where, for $i=\pm 1$, $\rA_i$ is a  copy of the vector space $\rA$ under the linear bijection
$x \mapsto x_i$.  Then
$\rL(\rA)$ is a Lie algebra,
called the
\emph{Tits-Kantor-Koecher Lie algebra}\index{Tits-Kantor-Koecher Lie algebra}
of $\rA$,
under the product defined by:
\begin{gather*}
[T,x_1] = (Tx)_1,\quad [T,y_{-1}] = (T^*y)_{-1},\quad  [x_i,y_i] = 0\ (i= \pm 1),\\
[x_1,y_{-1}] = V_{x,y},\quad  [T,T'] = TT' - T'T,
\end{gather*}
for $x,y\in\rA$, $T,T'\in  V_{\rA,\rA}$.
Finally, $\rL(\rA)$ is a $Q\times\Lm$-graded algebra
with
\[\rL(\rA)^\lm_{\al} = (\rA^\lm)_1, \quad
\rL(\rA)^\lm_0 = \sum_{\mu+\nu = \lm} V_{\rA^\mu,\rA^\nu}, \quad
\rL(\rA)^\lm_{-\al} = (\rA^\lm)_{-1},
\]
for $\lm\in\Lm$ and  $\rL(\rA)^\lm_{\beta} = 0$ for $\lm\in\Lm$, $\beta\in Q\setminus \De$.
Yoshii has shown  that $\rL(\rA)$ is a centreless Lie $\Lm$-torus of type
$\De$, and he proved the following coordinatization theorem.

\begin{theorem}[{\cite[Thm. 1]{Y1}}]
\label{thm:A1coord}
Any centreless Lie $\Lm$-torus of type $\De$
is graded-isomorphic to the Lie $\Lm$-torus  $\rL(\rA) = \TKK(\rA)$ constructed from a
Jordan $\Lm$-torus~$\rA$.
\end{theorem}

Yoshii  went on to  describe five families of Jordan  $\Zn$-tori and then
he showed that every
Jordan $\Lm$-torus\index{Jordan torus}
with $\rank(\Lm) = n$ is
isograded-isomorphic to a torus in one of the families \cite[Thm. 2]{Y1}.
The simplest of these families consists of the tori $\base_\bq^+$, where
$\base_\bq^+$ denotes the $\Zn$-graded algebra $\base_\bq$ with multiplication
$x\cdot y = \frac 12 (xy+yx)$.

\begin{definition}[Isotopes of Jordan tori]
\label{def:A1isotope}
Let  $\rA$ be Jordan $\Lm$-torus
and let $u$ be a nonzero homogeneous element of $\rA$.  So
$u\in \rA^{-\rho}$, where $\rho\in \supp_\lm(\rA)$.  Let
$\rA^\ui$ be the algebra with underlying vector space $\rA$ and product $\cdot_u$
defined by
\[x\cdot_u y = \frac 12 \{x,u,y\}.\]
Then, since $u$ is invertible in $\rA$, $\rA^\ui$ is a Jordan algebra and the identity
element of $\rA^\ui$ is $1^\ui = u^{-1}$ \cite[p.86]{Mc2}. So $1^\ui\in\rA^\rho$.
We endow $\rA^\ui$ with a
$\Lm$-grading by setting
\[(\rA^\ui)^\lm = \rA^{\lm + \rho}\]
for $\lm\in\Lm$. It is easily checked that $\rA^\ui$ is a Jordan $\Lm$-torus which
we call the
\emph{$u$-isotope}
\index{Jordan torus!isotope}
\index{Jordan torus!isotope!$\rA^\ui$}
of the Jordan torus $\rA$.  It is
also easily checked that up to graded-isomorphism $\rA^\ui$ does not depend
on the choice of nonzero $u$ in $\rA^{-\rho}$.
\end{definition}

Suppose that $u,v$ are nonzero homogeneous elements of a Jordan $\Lm$-torus.  It is well
known that
\begin{equation}
\label{eq:Jordisotope}
\rA^{(1)} = \rA,\quad (\rA^\ui)^{(v)} = \rA^{(U_uv)} \andd (\rA^\ui)^{(u^{-2})} = \rA,
\end{equation}
as algebras \cite[Proposition 7.2.1]{Mc2}, where $U_u v = \frac 12 \{u,v,u\}$.
It is easy to check that these
are also equalities of $\Lm$-graded Jordan algebras.

\begin{definition}
\label{def:Jordisotopy}
Suppose that $\rA$ is a Jordan $\Lm$-torus and $\rA'$ is a Jordan $\Lm'$-torus.
We say that $\rA$ is
\emph{isotopic}\index{Jordan torus!isotopic}
to  $\rA'$ if some
isotope of $\rA$ is isograded-isomorphic to $\rA'$.  In that case,
we write $\rA \sim \rA'$. It follows from \eqref{eq:Jordisotope} that
$\sim$ is an equivalence relation on the class of Jordan tori.
\end{definition}

The   following lemma is our key to understanding the connection between
isotopy for Lie tori and isotopy for Jordan tori.
In fact its proof (see \eqref{eq:A1discovery1} and
\eqref{eq:A1discovery2} below)
shows that the definition of isotope for Jordan tori
is determined without foreknowledge from the definition of isotope for
Lie tori.  We will see the same phenomenon for
alternative tori in \S \ref{sec:A2} and for associative tori with involution in~\S\ref{sec:Cr}.

\begin{lemma}[Key lemma for type A$_1$]
\label{lem:A1key}
Suppose that $\rA'$ is a Jordan $\Lm'$-torus, $\rA$ is
an Jordan $\Lm$-torus, and $\shom\in \Hom(Q,\Lm)$ is admissible for $\rL(\rA)$.
Suppose that $\ph : \rL(\rA') \to \rL(\rA)^\si$ is a bi-isomorphism with $\ph_r = 1$.
Then, $-\shom(\al) \in \supp_\Lm(\rA)$ and for any nonzero $u$ in $\rA^{-\shom(\al)}$
there is an isograded-isomorphism $\eta : \rA' \to \rA^\ui$
such that $\eta_\grtwist = \ph_e$.
\end{lemma}

\begin{proof} Note first that since $\shom$ is admissible, we have
$\shom(\al)\in\Lm_{\al}(\rL(\rA)) = \supp_\Lm(\rA)$, so $-\shom(\al)\in \supp_\Lm(\rA)$.
Thus we can choose nonzero $u$ in $\rA^{-\shom(\al)}$ as suggested.

Let $1'$ denote the identity in $\rA'$.  Then
\[\ph(1'_{-1}) \in \ph(\rL(\rA')_{-\al}^0) =
((\rL(\rA)^\si)_{-\al}^0 = \rL(\rA)_{-\al}^{-\shom(\al)} = (\rA^{-\shom(\al)})_{-1}.\]
So replacing $u$ by a scalar multiple we can assume that $\ph(1'_{-1}) = u_{-1}$.
Now define a linear bijection  $\eta : \rA' \to \rA$ by
\[(\eta x)_1 = \ph(x_1)\]
for $x\in\rA'$.
Then, for $x,y\in\rA'$, we have
\begin{align}
\notag
(\eta(xy))_1 &= \ph((xy)_1) =
\frac 12 \ph(\{x,1',y\}) =
\textstyle\frac 12 \ph([[x_1,1'_{-1}],y_1])
\\
\label{eq:A1discovery1}
&= \frac 12 [[(\eta x)_1,u_{-1}], (\eta y)_1]
= (\textstyle\frac 12 \{\eta x,u,\eta y\})_1 = (\eta(x)\cdot_u \eta(y))_1,
\end{align}
and, for $\lm\in\Lm'$,
\begin{align}
 \notag
(\eta({\rA'}^\lm))_1 &= \ph(({\rA'}^\lm)_1) = \ph(\rL(\rA')_{\al}^\lm)
 = (\rLs)_{\al}^{\ph_e(\lm)}
\\
\label{eq:A1discovery2}
&= (\rA^{\ph_e(\lm)+\shom(\al)})_1
 = ((\rA^\ui)^{\ph_e(\lm)})_1.
\end{align}\qed
\end{proof}

\begin{proposition}
\label{prop:A1isotope}  Let $\rA$ be a Jordan $\Lm$-torus,
and suppose that $\shom\in \Hom(Q,\Lm)$ is admissible for $\rL(\rA)$. Then,
$-\shom(\al) \in \supp_\Lm(\rA)$ and for any nonzero $u$ in $\rA^{-\shom(\al)}$, we have
$ \rL(\rA)^\si \simeq_{Q\times\Lm}\rL(\rA^\ui) $.
\end{proposition}

\begin{proof}  By  Theorem  \ref{thm:A1coord},  we have
$\rL(\rA)^\si \simeq_{Q\times\Lm} \rL(\rA')$ for some
Jordan $\Lm$-torus $\rA'$.  So we have a
bi-isomorphism $\ph : \rL(\rA') \to \rL(\rA)^\si$
with $\ph_r = 1$ and $\ph_e = 1$.
Lemma \ref{lem:A1key} now tells us that
$\rA' \simeq_\Lm \rA^\ui$,
so
$\rL(\rA')  \simeq_{Q\times \Lm} \rL(\rA^\ui)$.
Hence   $\rL(\rA)^\si  \simeq_{Q\times \Lm} \rL(\rA^\ui)$.
\qed\end{proof}

Our main theorem in this section, gives necessary and sufficient conditions for two
centreless Lie tori $\rL(\rA)$ and $\rL(\rA')$
to  be isotopic.  To do this, we need to first give necessary and sufficient
conditions for these Lie tori  to be bi-isomorphic.

\begin{theorem}
\label{thm:A1isotopy}
Suppose that $\rA$ is a Jordan $\Lm$-torus and $\rA'$ is a Jordan $\Lm'$-torus.
Let $\rL(\rA) = \TKK(\rA)$ and $\rL(\rA') = \TKK(\rA')$.
Then
\begin{description}[(ii) ]
\item[(i)] $\rL(\rA)\bi \rL(\rA') \iff \rA\ig\rA'$.
\item[(ii)] $\rL(\rA) \sim \rL(\rA') \iff \rA\sim\rA'$.
\end{description}
\end{theorem}

\begin{proof}  Let $\rL = \rL(\rA)$ and $\rL' = \rL(\rA')$.

(i): The proof of ``$\Leftarrow$'' is clear, so we prove ``$\Rightarrow$''.
Suppose $\ph : \rL' \to \rL$
is a bi-isomorphism.  Then $\ph_r = \pm 1$, so $\ph_r$ is in the Weyl group
of $\De$.
Therefore,  by Lemma \ref{lem:bi-isomorphism}, we can  assume that $\ph_r = 1$.
Thus, by Lemma \ref{lem:A1key} with $\shom = 0$ and $u=1$, we have
$\rA' \ig \rA$.

(ii): ``$\Rightarrow$''  By assumption, $\rL' \bi \rLs$ for some $\shom\in \Hom(Q,\Lm)$ that
is admissible for $\rL$.  So, choosing $u$ as in Proposition \ref{prop:A1isotope},
we have $\rL' \bi \rL(\rA^\ui)$.  By (i), $\rA' \ig \rA^\ui$.

``$\Leftarrow$'' By (i) we can assume
that $\rA' = \rA^\ui$, where $u$ is a nonzero homogeneous element of $\rA$.
Then $u\in\rA^{-\rho}$, where $\rho\in \supp_\Lm{\rL} = \Lm_{\al}(\rL)$.
So $\shom:Q\to \Lm$ defined by $\shom(\al) = \rho$ is admissible for $\rL$,
and hence $\rL' \bi \rL^\si$  by Proposition
\ref{prop:A1isotope}.
\qed\end{proof}

We will next give an example of two   centreless Lie tori that are isotopic but not bi-isomorphic.
For this we need a little preparation.

\begin{remark}
\label{rem:supportJordantori}  Suppose that $\rA$ is a Jordan $\Lm$-torus, and let
$\Ct(\rA)$ be the
centroid of $\rA$.  It is clear that
$\Ct(\rA)$ is $\Lm$-graded and that the support $\Gm = \Gm(\rA)$ of
$\Ct(\rA)$ is a subgroup of $\Lm$ \cite[\S 3]{Y1}.  Also, if
$S = \supp_\Lm(\rA)$, it is clear that $S$ is a union of cosets of $\Gm$ in
$\Lm$, and we let $S/\Gm = \set{\lm+\Gm \suchthat \lm\in S}\subseteq \Lm/\Gm$.
If $S/\Gm$ is finite,  we  denote the sum of the elements
of $S/\Gm$ in $\Lm/\Gm$ by $\Sigma(S/\Gm)$.
\end{remark}

\begin{lemma}
\label{lem:A1example}
Suppose $\rA$ is a Jordan $\Lm$-torus
with $S/\Gm$ finite and $\Sigma(S/\Gm) = 0$, where
$\Gm = \Gm(\rA)$ and $S = \supp_\Lm(\rA)$.
If $\rho\in S$ with  $\order{S/\Gm}\rho \notin \Gm$ and
$0 \ne u \in \rA^{-\rho}$, then $\rA^\ui \not\ig \rA$.
\end{lemma}

\begin{proof}   Let $\Gm^\ui = \Gm(\rA^\ui)$ and $S^\ui = \supp_\Lm(\rA^\ui)$.
It is easily checked that $\Ct(\rA^\ui) = \Ct(\rA)$,
so $\Gm^\ui = \Gm$.   Also, $S^\ui = S - \rho$, so
$\Sigma(S^\ui/\Gm) = \Sigma(S/\Gm) - \order{S/\Gm}(\rho + \Gm) = - \order{S/\Gm}(\rho + \Gm) \ne 0$.
But then $\Sigma(S/\Gm) = 0$ and $\Sigma(S^\ui/\Gm) \ne 0$ imply our conclusion.
\qed\end{proof}

Our example  uses Jordan tori of Clifford type \cite[Example 5.2]{Y1} as coordinate algebras.

\begin{example}
\label{ex:A1}
Let $\Lm = \bbZ^3$ and let $\rR = \base[2\Lm]$ be the  group algebra of $2\Lm$
with its  natural grading by $2\Lm$ (and hence by $\Lm$).  Let
$\set{\lm_1,\lm_2,\lm_3}$ be the standard basis of $\Lm$ and let
$\lm_4 = \lm_1+\lm_2+\lm_3$.  Let $\rV = \bigoplus_{i=1}^4 \rR x_i$
be the free $\rR$-module with base $\set{x_1,x_2,x_3,x_4}$, and give
$\rV$ the $\Lm$-grading such that $\deg(t^\mu x_i) = \mu+\lm_i$ for $\mu\in 2\Lm$ and
$1\le i \le 4$.  Define a $\Lm$-graded $\rR$-bilinear form $f:\rV \times \rV \to \rR$
by $f(x_i,x_j) = \delta_{i,j} t^{2\lm_i}$ for $1\le i,j\le 4$.  Finally, let $\rA = \rR\oplus \rV$
with the direct sum $\Lm$-grading, and define a product in $\rA$ by
\[(r+v)(r'+v') = rr' + f(v,v') + rv' + r'v\]
for $r,r'\in\rR$, $v,v'\in\rV$.
Then $\rA$ is a Jordan torus.  (The algebra
$\rA$ is called a \emph{spin factor} in \cite[p.74]{Mc2}.)
Further, using the notation of Lemma \ref{lem:A1example}, $\Gm = 2\Lm$ and
$S = (2\Lm)\cup (\cup_{i=1}^4 (2\Lm + \lm_i))$, so
$S/\Gm$ is finite  and $\Sigma(S/\Gm) = 0$.
Also
$\order{S/\Gm}\rho = 5\rho \notin \Gm$ for $\rho\in S\setminus \Gm$. So for any nonzero
homogeneous element $u$ in $\rA\setminus\rR$, we have, by Lemma \ref{lem:A1example}, that
 $\rA \not \ig \rA^\ui$. Therefore,
for any such $u$, we see by Theorem \ref{thm:A1isotopy} that
$\rL(\rA) \not\bi \rL(\rA^\ui)$ but $\rL(\rA) \sim \rL(\rA^\ui)$.
\end{example}

\section{Type A$_2$}
\label{sec:A2}

Suppose  in this section that
\[\De = \set{0}\cup \set{\ep_i-\ep_j \suchthat 1\le i  \ne j \le 3}\]
is the root system of type A$_2$ (where $\ep_1,\ep_2,\ep_3$
is a basis for a space containing $\De$), and $Q = Q(\De)$.  We let
$\al_1 = \ep_1-\ep_2$ and $\al_2 = \ep_2-\ep_3$, in which case
$\set{\al_1,\al_2}$ is a base for the root system $\De$.

Let $\rA$ be an alternative $\Lm$-torus.
We construct a centreless Lie $\Lm$-torus of type $\De$ from $\rA$
using J.~Faulkner's $3\times 3$-projective elementary
construction.  (See \cite[Appendix, (A6)]{F1} for the properties mentioned below about this construction.
Note however that we are using the notation $\eu_{ij}(x)$ in place of the notation $v_{ij}(x)$ that is used in
\cite{F1}.)
Recall that the
\emph{$3\times 3$ projective elementary Lie algebra}\index{projective elementary Lie algebra}
is the Lie algebra
\begin{equation}
\label{eq:pe3decomp}
\rL(\rA) = \pe_3(\rA) := \rL_0 \oplus \textstyle (\bigoplus_{1\le i \ne j \le 3} \eu_{ij}(\rA)),
\end{equation}
\index{projective elementary Lie algebra!$\pe_3(\rA)$}
where $\rL_0$ and $\eu_{ij}(\rA)$ are subspaces of $\rL$ satisfying:
\begin{equation}
\label{eq:pe3prop}
\parbox[l]{4.5in}
{
\begin{description}[(a) ]
\item[(a)] For $i\ne j$, there is a linear bijection  $x \mapsto \eu_{ij}(x)$ from $\rA$ to $\eu_{ij}(\rA)$,
\item[(b)] $[\eu_{ij}(x),\eu_{jk}(y)] = \eu_{ik}(xy)$ for $x,y\in \rA$ and distinct $i,j,k$,
\item[(c)] $[\eu_{ij}(x),\eu_{ij}(y)] = 0$  for  $x,y\in \rA$, $i\ne j$,
\item[(d)] $[\rL_0,\rL_0] \subseteq \rL_0$,
\item[(e)] $[\rL_0,\eu_{ij}(\rA)] \subseteq \eu_{ij}(\rA)$ for $i\ne j$,
\item[(f)]$\rL_0 = \sum_{1\le i< j \le 3} [\eu_{ij}(\rA),\eu_{ji}(\rA)]$,
\item[(g)] $\set{x\in \rL_0: [x,\eu_{ij}(\rA)] = 0 \text{ for } i\ne j} = 0$.
\end{description}
}
\end{equation}
These properties characterize the Lie algebra $\rL(\rA)$ in the following sense:  If
$\rL'(\rA) = \rL'_0  \oplus \allowbreak \textstyle (\bigoplus_{1\le i \ne j \le 3} \eu'_{ij}(\rA))$
is another Lie algebra satisfying the conditions in \eqref{eq:pe3prop}, then there
is a unique isomorphism from $\rL(\rA) \to \rL'(\rA)$ such that
$\eu_{ij}(x) \mapsto \eu'_{ij}(x)$ for $i\ne j$, $x\in \rA$.
Finally $\rL(\rA)$ is a $Q\times \Lm$-graded algebra with
\[\rL(\rA)_{\ep_i-\ep_j}^\lm = \eu_{ij}(\rA^\lm),\quad
\rL(\rA)_0^\lm = \textstyle \sum_{\mu+\nu = \lm} \sum_{1\le k< l \le 3} [\eu_{kl}(\rA^\mu),\eu_{lk}(\rA^\nu)],\]
for $i\ne j$, $\lm\in \Lm$, and $\rL(\rA)_\al^\lm = 0$ for $\lm\in \Lm$, $\al\in Q\setminus \De$.
It is easily checked that $\rL(\rA)$ is a centreless Lie torus of type A$_2$.

Moreover, Berman, Gao, Krylyuk and  Neher (over $\bbC$) and Yoshii (in general)
proved the following coordinatization theorem:

\begin{theorem} \emph{(\cite[Lemma 3.25]{BGKN}, \cite[Prop.~6.3]{Y4})}
\label{thm:A2coord}
Any  centreless Lie $\Lm$-torus of type $\De$
is grad\-ed-isomorphic to the Lie $\Lm$-torus $\rL(\rA) = \pe_3(\rA)$ constructed from an
alternative $\Lm$-torus $\rA$.\footnote{A different description of
$\rL(\rA)$ was used in \cite{BGKN} and \cite{Y4}, but it is not difficult to check that
the two descriptions give Lie tori that
are $Q\times \Lm$-graded-isomorphic.}
\end{theorem}

An important example of an alternative torus is the octonion torus.  A simple way
to introduce this torus uses a presentation \cite{Y7}.

\begin{example}  Let $\cC(3)$ be the alternative algebra  that is presented by the generators
$x_i^{\pm 1}, 1\le i \le 3$,  subject to the relations $x_ix_i^{-1} = x_i^{-1}x_i = 1$ for all $i$,
\[x_ix_j = -x_j x_i \text{ for } i\ne j,\andd (x_1x_2)x_3 = -x_1(x_2x_3).\]
Although we will not need this, one can show more concretely that each element
of $\cC(3)$ can be written uniquely as a linear combination of monomials
$x^{\boldsymbol i} = (x_1^{i_1}x_2^{i_2})x_3^{i_3}$, $\boldsymbol i = (i_1,i_2,i_3)\in \bbZ^3$,
with multiplication given by
\[x^{\boldsymbol i} x^{\boldsymbol j}
= (-1)^{\kappa(\boldsymbol i,\boldsymbol j)}
x^{\boldsymbol i + \boldsymbol j},
\]
where $\kappa(\boldsymbol i,\boldsymbol j) = i_3j_1 + i_2j_1+i_3j_2 + i_1j_2j_3 + i_2j_1j_3+ i_3j_1j_2$.

If $n\ge 3$, the graded algebra $\cC(3) \otimes \base[x_4^{\pm 1},\dots,x_n^{\pm 1}]$ is
an alternative $\bbZ^n$-torus, where
$x_1,\dots,x_n$ are assigned the degrees
in the standard ordered basis for $\Zn$ (identifying
$\cC(3)$ and $\base[x_4^{\pm 1},\dots,x_n^{\pm 1}]$ as subalgebras
of the tensor product).  The $\Zn$-graded algebra
$\cC(3) \otimes \base[x_4^{\pm 1},\dots,x_n^{\pm 1}]$ is called
the
\emph{octonion $\bbZ^n$-torus}\index{octonion torus}.
\end{example}

We  have the following description of
alternative tori\index{alternative torus}:

\begin{theorem} \emph{(\cite[Thm.~1.25]{BGKN}, \cite[Cor.~5.13]{Y4})}
\label{thm:A2class}
If $\rA$ is an  alternative $\Lm$-torus,
then either $\rA$ is associative (and hence isograded-isomorphic to a quantum torus),
or $n\ge 3$ and $\rA$ is isograded-isomorphic to the octonion $\bbZ^n$-torus.
\end{theorem}

\begin{remark}
\label{rem:A2support}  Suppose $\rA$ is a $\Lm$-torus.
Recall that by
Remark \ref{rem:supporttori}, $\supp(\rA) = \Lm$.
So $\Lm_\al(\rL(\rA)) = \Lm$ for all nonzero $\al\in\De$, and therefore
any $\shom\in \Hom(Q,\Lm)$ is admissible for $\rL(\rA)$.
\end{remark}

\begin{definition}[Isotopes of alternative tori]
\label{def:A2isotope}
Let $\rA$ be
an alternative $\Lm$-torus and let $u_1$ and $u_2$ be nonzero homogeneous elements
of $\rA$.  Thus, $u_i\in \rA^{-\rho_i}$, where $\rho_i \in \Lm$ for $i=1,2$.
Let $\rA^{(u_1,u_2)}$
be the algebra with underlying vector space $\rA$ and with product $\cdot_{u_1,u_2}$
defined by:
\[ x \cdot_{u_1,u_2} y = (xu_1)(u_2 y)\]
for $x,y\in\rA$.  Then, since $u_1,u_2$ are invertible,  $\rA^{(u_1,u_2)}$ is an alternative
algebra with identity $1^{(u_1,u_2)} = (u_1 u_2)^{-1}$ \cite[Thm.~1 and Prop.~2]{Mc1}.  We define a $\Lm$-grading
on $\rA^{(u_1,u_2)}$ by setting
\[(\rA^{(u_1,u_2)})^\lm = \rA^{\lm + \rho_1 + \rho_2}\]
for  $\lm\in \Lm$.  Then $\rA^{(u_1,u_2)}$ is an alternative $\Lm$-torus which we call
the
\emph{$(u_1,u_2)$-isotope of $\rA$}\index{alternative torus!isotope}
\index{alternative torus!isotope!$\rA^{(u_1,u_2)}$}.
It is easily checked that this torus does not depend
up to graded-isomorphism on the choice of nonzero $u_i$ in $\rA^{-\rho_i}$.
\end{definition}

\begin{lemma}[Key lemma for type A$_2$]
\label{lem:A2key}
Suppose that $\rA'$ is an alternative $\Lm'$-torus, $\rA$ is
an alternative $\Lm$-torus, and $\shom\in \Hom(Q,\Lm)$.
Suppose that $\ph : \rL(\rA') \to \rL(\rA)^\si$ is a bi-isomorphism with $\ph_r = 1$.
Choose $0\ne u_i\in \rA^{\shom(\al_i)}$
for $i = 1,2$. Then there is an isograded-isomorphism $\eta : \rA' \to \rA^{(u_1,u_2)}$
such that $\eta_\grtwist = \ph_e$.
\end{lemma}

\begin{proof}  Write $\rL(\rA') = \rL'_0 \oplus (\bigoplus_{1\le i \ne j \le 3} \eu'_{ij}(\rA'))$
as in \eqref{eq:pe3decomp}.
Now $\ph(\eu'_{ij}(\rA'^\lm)) = \ph(\rL(\rA')_{\ep_i-\ep_j}^\lm)
= (\rL(\rA)^\si)_{\ep_i-\ep_j}^{\ph_e(\lm)} = \rL(\rA)_{\ep_i-\ep_j}^{\ph_e(\lm)+ \shom(\ep_i-\ep_j)} =  \eu_{ij}(\rA^{\ph_e(\lm) + \shom(\ep_i-\ep_j)})$
for $\lm\in\Lm'$, $i\ne j$.  Thus, for $i\ne j$,
we have a linear bijection $\ph_{ij} : \rA' \to \rA$ with
\begin{equation}
\label{eq:A2phij1}
\ph(\eu'_{ij}(x)) = \eu_{ij}(\ph_{ij}(x))
\end{equation}
for $x\in\rA'$, in which case
\begin{equation}
\label{eq:A2phij2}
\ph_{ij}({\rA'}^\lm) =  \rA^{\ph_e(\lm)+ \shom(\ep_i-\ep_j)}
\end{equation}
for $\lm\in\Lm'$.
Using \eqref{eq:A2phij1} and the relation \eqref{eq:pe3prop}(b) in $\rL(\rA')$ and $\rL(\rA)$,
we see that
\begin{equation}
\label{eq:A2phij3}
\ph_{ij}(xy) = \ph_{ik}(x)\ph_{kj}(y)
\end{equation}
for $x,y\in\rA'$ and $i,j,k\ne$.

Now let $u_{ij} = \ph_{ij}(1')$ for $i\ne j$, where $1'$ is the identity in $\rA'$.
Then, taking $y=1'$ and separately $x=1'$ in \eqref{eq:A2phij3}, we have $\ph_{ij}(x) = \ph_{ik}(x) u_{kj}$ and $\ph_{ij}(y) = u_{ik}\ph_{kj}(y)$
for $x,y\in\rA'$ and $i,j,k\ne$.  So
\[\ph_{ij}(xy) = \ph_{ik}(x)\ph_{kj}(y) = (\ph_{ij}(x)u_{jk})(u_{ki}\ph_{ij}(y))
= \ph_{ij}(x)\cdot_{u_{jk},u_{ki}}\ph_{ij}(y)
\]
for $x,y\in\rA'$ and $i,j,k\ne$.  Also, by \eqref{eq:A2phij2}, $u_{12}\in \rA^{\shom(\al_1)}$
and $u_{23}\in \rA^{\shom(\al_2)}$.  So replacing $u_1$ and $u_2$ by scalar multiples,
we may assume that $u_1 = u_{12}$ and $u_2 = u_{23}$.  Thus, $\ph_{31}$ is an algebra isomorphism
of $\rA'$ onto $\rA^{(u_1,u_2)}$.  Moreover, for $\lm\in\Lm'$,
\[\ph_{31}((\rA')^\lm) = \rA^{\ph_e(\lm) + \shom(\ep_3-\ep_1)} = \rA^{\ph_e(\lm) - \shom(\al_1)-\shom(\al_2)}
= (\rA^{(u_1,u_2)})^{\ph_e(\lm)}.\]
So we have our conclusion with $\eta = \ph_{31}$.
\qed\end{proof}

\begin{proposition}
\label{prop:A2isotope}  Suppose  $\rA$ is  an alternative $\Lm$-torus
and $\shom\in \Hom(Q,\Lm)$. Choose $0\ne u_i\in \rA^{\shom(\al_i)}$
for $i = 1,2$. Then,
$ \rL(\rA)^\si \simeq_{Q\times\Lm}\rL(\rA^{(u_1,u_2)})$.
\end{proposition}

\begin{proof}  This follows from Lemma \ref{lem:A2key} as in \S \ref{sec:A1}.
(See the proof of Proposition \ref{lem:A1key}.)
\qed\end{proof}

We could introduce isotopy for
alternative tori\index{alternative torus!isotopic}
in the obvious fashion (as in Definition
\ref{def:Jordisotopy} in the Jordan case).
However, in view of the following result, this will not be necessary.

\begin{proposition}
\label{prop:A2isotopytrivial}
Suppose that $\rA$ is an alternative $\Lm$-torus and $u_1,u_2$ are nonzero homogeneous
elements of $\rA$. Then $\rA^{(u_1,u_2)} \ig \rA$.
\end{proposition}

\begin{proof} Let $\rA' = \rA^{(u_1,u_2)}$.
Suppose first that $\rA$ is associative.  Then, one checks that the map
$x \mapsto (u_1u_2)^{-1}x$ is a graded algebra isomorphism of $\rA$ onto $\rA'$.
So we can assume that $\rA$ is not associative.  Thus,
$\rA'$ is also not associative \cite[p.256--257]{Mc1}.
Therefore, by Theorem \ref{thm:A2class},  $\rA$ and $\rA'$
are both isograded-isomorphic to the octonion $\Zn$-torus.  So $\rA \ig \rA'$.
\qed\end{proof}

\begin{remark} With the assumptions of Proposition \ref{prop:A2isotopytrivial} one can actually show
that $\rA^{(u_1,u_2)} \simeq_\Lm \rA$, but isograded isomorphism is all that we need.
\end{remark}

If $\rA$ is an alternative $\Lm$-torus, then the
\emph{opposite torus}\index{alternative torus!opposite}
of $\rA$
is the alternative $\Lm$-torus $\rA^\opalg$, where $\rA^\opalg = \rA$ as graded vector spaces and the product on $\rA^\opalg$ is given by $x\cdot_\opalg y = yx$.

\begin{lemma}
\label{lem:A2opposite}  Suppose that $\rA$ is an alternative $\Lm$-torus.  Then, there
exists a unique bi-isomorphism from $\rL(\rA)$ onto $\rL(\rA^\opalg)$ such that
\[\eu_{ij}(x) \mapsto -\tilde \eu_{ji}(x)\]
for $x\in \rA$ and $i\ne j$,
where $\rL(\rA^\opalg) = \tilde\rL_0 \oplus (\bigoplus_{1\le i \ne j \le 3} \tilde \eu_{ij}(\rA^\opalg))$
as in \eqref{eq:pe3decomp}.
\end{lemma}

\begin{proof}  In $\rL(\rA^\opalg)$, put
$\eu'_{ij}(x) = -\tilde \eu_{ji}(x)$ for $x\in \rA$, $i\ne j$,
and put $\rL'_0 = \tilde\rL_0$.
Then, $\rL(\rA^\opalg) = \rL'_0 \oplus (\bigoplus_{1\le i \ne j \le 3} \eu'_{ij}(\rA^\opalg))$,
and for $x,y\in\rA$ and $i,j,k\ne$, we have
\[[\eu'_{ij}(x),\eu'_{jk}(y)] = [\tilde \eu_{ji}(x),\tilde \eu_{kj}(y)] =
-[\tilde \eu_{kj}(y),\tilde \eu_{ji}(x)] = -\tilde \eu_{ki}(y\cdot_\opalg x) = \eu'_{ik}(xy).
\]
The other properties in \eqref{eq:pe3prop} are clear, so
we have a unique Lie algebra isomorphism
$\ph : \rL(\rA) \to \rL(\rA^\opalg)$ such that
$\ph(\eu_{ij}(x)) = \eu'_{ij}(x) =  -\tilde \eu_{ji}(x)$ for $x\in \rA$ and $i\ne j$.  It is clear
that $\ph$ is a bi-isomorphism with $\ph_r = -1$ and $\ph_e = 1$.
\qed\end{proof}

The main result of this section is:

\begin{theorem}
\label{thm:A2isotopy}
Suppose that $\rA$ is an alternative $\Lm$-torus and $\rA'$ is an alternative $\Lm'$-torus.
Let $\rL(\rA) = \pe_3(\rA)$ and $\rL(\rA') = \pe_3(\rA')$. Then,
the following are
equivalent:\footnote{See \cite[Lemma 1.4]{GN} for a result related to Theorems \ref{thm:A2isotopy}
and \ref{thm:Arisotopy} in the context of Jordan pairs covered by grids.}
\begin{description}[(a) ]
\item[(a)] $\rA$ is isograded-isomorphic to $\rA'$ or $(\rA')^\opalg$.
\item[(b)] $\rL(\rA) \bi \rL(\rA')$.
\item[(c)] $\rL(\rA)\sim \rL(\rA')$.
\end{description}
\end{theorem}

\begin{proof}  ``(a)$\Rightarrow$(b)''  Since $\rL(\rA') \bi \rL({\rA'}^\opalg)$ by
Lemma \ref{lem:A2opposite}, we can assume that $\rA \ig \rA'$.  (b) is then clear.

``(b)$\Rightarrow$(c)'' is trivial.

``(c)$\Rightarrow$(a)'' By assumption, $\rL(\rA') \bi \rL(\rA)^\si$ for some
$\shom\in \Hom(Q,\Lm)$ that is admissible for $\rL(\rA)$.
Choose $0\ne u_i\in \rA^{\shom(\al_i)}$ for $i=1,2$.
Then, by Propositions \ref{prop:A2isotope} and \ref{prop:A2isotopytrivial},
we have $\rL(\rA)^\si \bi \rL(\rA^{(u_1,u_2)}) \bi \rL(\rA)$.
So there is a bi-isomorphism $\ph: \rL(\rA') \to \rL(\rA)$.
Now
the automorphism group of $\De$ is generated by $-1$ and the Weyl group of $\De$.  Hence,
post-multiplying $\ph$ by the bi-isomorphism in Lemma \ref{lem:A2opposite}
if necessary, we can assume that $\ph_r$ is in the Weyl group.
Thus, by Lemma \ref{lem:bi-isomorphism}, we can assume that $\ph_r = 1$.
So, by Lemma \ref{lem:A2key} (with $\shom = 0$ and $u_1=u_2=1$), we have $\rA' \ig \rA$.
\qed\end{proof}

\begin{remark}
\label{rem:opposite}
If $\rA$ is an alternative torus, it is not in general true that $\rA$ is isograded-isomorphic
to $\rA^\opalg$.  However, one can show that this is true if $\rA$ is the octonion torus
(since the octonion torus has an involution) and if $\rA$ is a rational quantum torus
(using the classification in   \cite{Ne}).
\end{remark}

\section{Type A$_r$, $r\ge 3$}
\label{sec:Ar}

Suppose in this section that $r\ge 3$,
\[\De = \set{0}\cup \set{\ep_i-\ep_j \suchthat 1\le i  \ne j \le r+1}\]
is the root system of type A$_r$ (where $\ep_1,\dots,\ep_r$
is a basis for a space containing $\De$), and $Q = Q(\De)$.
We let $\al_i = \ep_i-\ep_{i+1}$ for $1\le i \le r$, in which case
$\set{\al_1,\dots,\al_r}$ is a base for the root system $\De$.

Let $\rA$ be an associative $\Lm$-torus, and denote by $[\rA,\rA]$
the Lie algebra spanned by commutators in $\rA$. Let
\[\rL(\rA) = \speciallinear_{r+1}(\rA) :=
\set{X\in \Mat_{r+1}(\rA) \suchthat \trace(X) \in [\rA,\rA] }.\]
\index{special linear Lie algebra!$\speciallinear_{r+1}(\rA)$}
Then $\rL(\rA)$ is a Lie algebra under the commutator product called
the
\emph{$(r+1)\times (r+1)$-special linear}
\index{special linear Lie algebra}
Lie algebra over $\rA$.
For $1\le i,j\le r+1$, we let $\eu_{ij}$ be the \emph{$(i,j)$-matrix unit} in $\Mat_{r+1}(\rA)$
with  $1$ in the $(i,j)$-position and zeros elsewhere,
and we let  $\eu_{ij}(x) = xe_{ij}$ for $x\in \rA$.
Then $\rL(\rA)$ is a $Q\times \Lm$-graded algebra with
\[\rL(\rA)_{\ep_i-\ep_j}^\lm = \eu_{ij}(\rA^\lm),\quad
\rL(\rA)_0^\lm =
\set{\textstyle\sum_{i=1}^{r+1} \eu_{ii}(x_i) \suchthat x_i \in \rA^\lm,\ \sum_{i=1}^{r+1} x_i \in [\rA,\rA]},\]
for $i\ne j$, $\lm\in \Lm$, and $\rL(\rA)_\al^\lm = 0$ for $\lm\in \Lm$, $\al\in Q\setminus \De$.
It is easy to check that the matrices $\eu_{ij}(x)$, $x\in\rA$, $i\ne j$, generate the Lie algebra
$\rL(\rA)$ and also, using the fact that the centre of $\rA$ intersects trivially
with $[\rA,\rA]$ \cite[Prop.~2.44]{BGK},
that $\rL(\rA)$ is
centreless.\footnote{In view of these properties, it is natural to also use the notation
$\pe_{r+1}(\rA)$ for $\rL(\rA)$, which is consistent with the notation in rank 2 (see Remark
\ref{rem:Arcompareranks}).}
It then follows easily that
$\rL(\rA)$ is a centreless Lie $\Lm$-torus of type $\De$.

Moreover, Berman, Gao and  Krylyuk   (over $\bbC$) and Yoshii (in general)
proved the following coordinatization theorem:

\begin{theorem} \emph{(\cite[Thm.~2.65]{BGK}, \cite[\S 4]{Y2})}
\label{thm:Arcoord}
Any centreless Lie $\Lm$-torus of type $\De$
is grad\-ed-iso\-mor\-phic to the Lie $\Lm$-torus $\rL(\rA)= \speciallinear_{r+1}(\rA)$ constructed from an
associative $\Lm$-torus $\rA$.
\end{theorem}

\begin{remark}
\label{rem:Arcompareranks}
Exactly as above one can define
centreless Lie tori $\speciallinear_2(\rA)$ and $\speciallinear_3(\rA)$ in rank 1 and 2 respectively.
Then it is not difficult to show that $\speciallinear_2(\rA) \simeq_{Q\times \Lm} \TKK(\rA^+)$ and
$\speciallinear_3(\rA) \simeq_{Q\times \Lm} \pe_3(\rA)$.  However, as we've seen,
these constructions do not provide all examples of Lie tori of type A$_1$ or A$_2$.
\end{remark}

To understand isotopes of  Lie tori of type $\De$, we could proceed as in rank 1 and 2 to
relate them to isotopes of the coordinate tori.  However, we can take a simpler approach
in view of the next  proposition.\footnote{This proposition
could also be proved at the level of coordinates using the fact that
any isotope of an associative algebra $\rA$ is isomorphic
to $\rA$.}

\begin{proposition}
\label{prop:Arisotopytrivial}
Suppose  $\rA$ is an associative $\Lm$-torus and
$\shom\in \Hom(Q,\Lm)$.  (Note that as in rank 2
any $\shom\in \Hom(Q,\Lm)$ is admissible for $\rL(\rA)$.)
Then,
$\rL(\rA)^\si \simeq_{Q\times\Lm} \rL(\rA)$.
\end{proposition}

\begin{proof} Let $Q' = \bigoplus_{i=1}^{r+1} \bbZ\ep_i$.  Then
$Q' = \bbZ\ep_1 \oplus Q$, so we can extend $\shom$ to a group
homomorphism $\shom: Q'\to \Lm$
by setting $\shom(\ep_1) = 0$.  Choose $0\ne d_i\in \rA^{\shom(\ep_i)}$
for $1\le i \le r+1$.  Let $d = \diag(d_1,\dots,d_{r+1})$
and define $\ph : \rL(\rA) \to \rL(\rA)$ by
\[\ph(X) = dXd^{-1}\]
for $X\in \rL(\rA)$.  (Here we are using the fact that
$\trace(X_1X_2) \equiv \trace(X_2X_1) \pmod{[\rA,\rA]}$ for $X_1,X_2\in\Mat_{r+1}(\rA)$,
so $\ph$ maps $\rL(\rA)$ into $\rL(\rA)$.)  Then, $\ph$ is a Lie algebra isomorphism.
Also
\begin{align*}
\ph(\rL(\rA)_{\ep_i-\ep_j}^\lm) &= \ph(\eu_{ij}(\rL^\lm)) = \eu_{ij}(d_i\rA^\lm d_j^{-1})
= \eu_{ij}(\rA^{\lm + \shom(\ep_i) - \shom(\ep_j)}) \\
&= \rL(\rA)_{\ep_i-\ep_j}^{\lm+\shom(\ep_i-\ep_j)} =
(\rL(\rA)^\si)_{\ep_i-\ep_j}^\lm
\end{align*}
for $\lm\in\Lm$; so, by (LT3), $\ph(\rL(\rA)_\al^\lm) = (\rL(\rA)^\si)_\al^\lm$ for $\al\in Q$, $\lm\in\Lm$.
\qed\end{proof}

The ``key lemma'' that we need is simpler  than in rank 1 and 2
because of  Proposition \ref{prop:Arisotopytrivial}.

\begin{lemma}[Key lemma for type A$_r$, $r\ge 3$]
\label{lem:Arkey}
Suppose that $\rA'$ is an associative $\Lm'$-torus and $\rA$ is
an associative $\Lm$-torus.
Suppose that $\ph : \rL(\rA') \to \rL(\rA)$ is a bi-isomorphism with $\ph_r = 1$.
Then, there is an isograded-isomorphism $\eta : \rA' \to \rA$
such that $\eta_\grtwist = \ph_e$.
\end{lemma}

\begin{proof} By Lemma \ref{lem:bi-isomorphism}, we can assume that
$\ph(\eu_{ij}(1')) = \eu_{ij}(1)$ for $i-j = \pm 1$, and hence for
all $i\ne j$.  Also, we have the relation
\[[\eu_{ij}(x),\eu_{jk}(y)] = \eu_{ik}(xy)\]
in $\rL(\rA)$ for $x,y\in\rA$ and $i,j,k\ne$.
We can thus argue as in the proof of Lemma
\ref{lem:A2key} (with $\shom = 0$ and $u_1=u_2=1$).
\qed\end{proof}

The following is easily checked:

\begin{lemma}
\label{lem:Aropposite}  Suppose that $\rA$ is an associative $\Lm$-torus.  Then,
the map $T \mapsto -T^t$ is a bi-isomorphism from $\rL(\rA)$ onto $\rL(\rA^\opalg)$.
\end{lemma}

\begin{theorem}
\label{thm:Arisotopy}
Suppose that $\rA$ is an associative $\Lm$-torus and $\rA'$ is an associative $\Lm'$-torus.
Let $\rL(\rA) = \speciallinear_{r+1}(\rA)$ and $\rL(\rA') = \speciallinear_{r+1}(\rA')$.
Then,
the following are equivalent:
\begin{description}[(a) ]
\item[(a)] $\rA$ is isograded-isomorphic to $\rA'$ or $(\rA')^\opalg$.
\item[(b)] $\rL(\rA) \bi \rL(\rA')$.
\item[(c)] $\rL(\rA)\sim \rL(\rA')$.
\end{description}
\end{theorem}

\begin{proof}  If (a) holds then, by Lemma \ref{lem:Aropposite}, we can assume that $\rA \ig \rA'$,
so (b) follows. The implication ``(b)$\Rightarrow$(c)'' is trivial. Finally, suppose that
(c) holds. Thus, $\rL(\rA)^\si \bi \rL(\rA')$ for some
$\shom\in \Hom(Q,\Lm)$, so
$\rL(\rA) \bi \rL(\rA')$ by Proposition \ref{prop:Arisotopytrivial}.
We can now obtain (a) using the argument in the proof of Theorem \ref{thm:A2isotopy}
(and using Lemmas \ref{lem:Arkey} and \ref{lem:Aropposite}).
\qed\end{proof}

\section{Type C$_r$, $r\ge 4$}
\label{sec:Cr}

Suppose in this section that $r\ge 4$,
\[\De = \set{0}\cup \set{\ep_i-\ep_j \suchthat 1\le i  \ne j \le r}
 \set{\pm(\ep_i+\ep_j) \suchthat 1\le i  \le j \le r}\]
is the root system of type C$_r$ (where $\ep_1,\dots,\ep_r$
is a basis for the space spanned by $\De$), and $Q = Q(\De)$.
We let $\al_i = \ep_i-\ep_{i+1}$, $1\le i \le r-1$, and $\al_r = 2\ep_r$,
in which case
$\set{\al_1,\dots,\al_r}$ is a base for the root system $\De$.

Let $(\rA,\iota)$ be an associative $\Lm$-torus with involution.  For convenience, write
$\bar x = \iota(x)$ for $x\in \rA$ and $\bar X = (\bar x_{ij})$ for $X = (x_{ij}) \in \Mat_r(\rA)$.
We write $\rA_\pm =  \rA_{\pm,\,\iota} := \set{x\in \rA \suchthat \bar x = \pm x}$, in which case $\rA = \rA_+ \oplus \rA_-$.
Let
\begin{align*}
\rL(\rA,\iota) = \ssp_{2r}(\rA,\iota) :=
\set{\begin{bmatrix} X&Y \\ Z & -\bar X^t \end{bmatrix} &\suchthat X,Y,Z\in \Mat_r(\rA),\ \\
&\bar Y^t = Y,\ \bar Z^t = Z,\ \trace(X) \in \rA_+ + [\rA,\rA] }.
\end{align*}
\index{special symplectic Lie algebra!$\ssp_{2r}(\rA,\iota)$}%
Then $\rL(\rA,\iota)$ is a Lie algebra under the commutator product called
the
\emph{$(2r)\times (2r)$-special symplectic Lie algebra}
\index{special symplectic Lie algebra}
over $\rA$.
To define a grading on $\rL(\rA,\iota)$, let
\begin{gather*}
\ell_{ij}(x) = x\eu_{i,j} - \bar x \eu_{j+r,i+r} \quad\text{ if } 1\le i \ne j \le r,\ \text{ and }\\
m_{ij}(x) = x\eu_{i,j+r} + \bar x \eu_{j,i+r},\quad  n_{ij}(x) = x\eu_{i+r,j} + \bar x \eu_{j+r,i}
\quad\text{ if } 1\le i, j \le r,
\end{gather*}
where  $\eu_{ij}$ is the $(i,j)$-matrix unit.
With this notation, $\rL(\rA,\iota)$ is a $Q\times \Lm$-graded algebra with
\begin{align*}
\rL(\rA,\iota)_{\ep_i-\ep_j}^\lm &= \ell_{ij}(\rA^\lm) \quad \text{for } 1 \le i \ne j \le r,\\
\rL(\rA,\iota)_{\ep_i+\ep_j}^\lm &= m_{ij}(\rA^\lm) \quad\text{for }  1 \le i,j \le r,\\
\rL(\rA,\iota)_{-(\ep_i+\ep_j)}^\lm &= n_{ij}(\rA^\lm) \quad \text{for } 1\le i,j \le r,
\end{align*}
\begin{align*}
\rL(\rA,\iota)_0^\lm &=
\set{\textstyle\sum_{i=1}^{r} (\eu_{i,i}(x_i)- \eu_{i+r,i+r}(\bar x_i)) \suchthat
x_i \in \rA^\lm,\ \sum_{i=1}^r x_i \in \rA_+ + [\rA,\rA]}
\end{align*}
for $\lm\in \Lm$, and $\rL(\rA,\iota)_\al^\lm = 0$ for $\lm\in \Lm$, $\al\in Q\setminus \De$.
Then $\rL(\rA,\iota)$ is a centreless Lie $\Lm$-torus of type~$\De$.

Moreover, Allison and Gao (over $\bbC$) and Benkart and  Yoshii (in general) proved:

\begin{theorem} \emph{(\cite[Prop.~4.87]{AG}, \cite[Theorem 5.9]{BY})}
\label{thm:Crcoord}
Any centreless Lie $\Lm$-torus of type $\De$
is grad\-ed-iso\-mor\-phic to the Lie $\Lm$-torus $\rL(\rA,\iota) = \ssp_{2r}(\rA,\iota)$ constructed from an
associative $\Lm$-torus
with involution  $(\rA,\iota)$.\footnote{$\rL(\rA,\iota)$ was described differently
in \cite{AG} and also in \cite{BY}.  It is not difficult to check that the three descriptions give Lie tori that
are $Q\times \Lm$-graded-isomorphic.}
\end{theorem}

\begin{remark}
\label{rem:Crcompareranks}
Exactly as above one can define
centreless Lie tori $\ssp_{4}(\rA,\iota)$ and  $\ssp_{6}(\rA,\iota)$.
However, these constructions do not provide all examples of Lie tori of type C$_2 = \text{B}_2$  or C$_3$
\cite[Thm.~4.87]{AG}.
\end{remark}

\begin{example}
\label{ex:quantuminvolution}
Let $\bq = (q_{ij})\in \Mat_n(\base)$ with
$q_{ii} = 1$ and $q_{ij} = q_{ji} = \pm 1$ for all $i,j$,
and let
$\base_\bq = \base_\bq[x_1^{\pm1},\dots,x_n^{\pm1}]$  be the quantum torus
corresponding to $\bq$ as in Example \ref{ex:quantum}.
Suppose $\be = (e_1,\dots,e_n)\in \base^n$ with $e_i = \pm 1$,   Then,
there is a unique involution $\iota_\be$ on $\base_\bq$ such that
$\iota_\be(x_i) = e_i x_i$ for $1\le i \le n$; and
$(\base_\bq,\iota_\be)$ is an associative $\Zn$-torus with
involution.\footnote{If $n=0$, we interpret
$(\base_\bq,\iota_\be)$ as $(\base,\id)$.}
It is not difficult to show that any associative torus with involution $(\rA,\iota)$
is isograded-isomorphic to some $(\base_\bq,\iota_\be)$.  More specifically,
if we choose any basis $\set{\lm_1,\dots,\lm_n}$ for $\Lm$ and choose  $0\ne y_i\in \rA^{\lm_i}$
for each $i$, then there is an isograded-isomorphism of $(\rA,\iota)$ onto
$(\base_\bq,\iota_\be)$ so that $y_i \mapsto x_i$ for each $i$, where
$\bq$ and $\be$ are determined by the equalities $y_jy_i = q_{ij} y_i y_j$ and
$\iota(y_i) = e_i y_i$ \cite[p.163--164]{AG}.
\end{example}

Yoshii in \cite[Thm. 2.7]{Y3}  classified the
associative tori with involution\index{associative torus with involution}
$(\base_\bq,\iota_\be)$ up to isograded-isomorphism using successive changes of variables
that transform the pair $(\bq,\be)$ into a canonical form.  Alternatively,  as we now see,
one can classify associative tori with involution using quadratic forms
following an approach discussed in \cite[Remark 5.20]{AFY}.

Let $\bLm = \Lm/2\Lm$ regarded as a vector space over $\Ztwo = \bbZ/2\bbZ$,
and let $\lm \to \blm$ be the canonical map of $\Lm$ onto $\bLm$.  Recall that
a map $\kappa : \bLm \to \Ztwo$ is called a \emph{quadratic form} on $\bLm$
if the map $\kappa_p : \bLm \times \bLm \to \Ztwo$ defined by
$\kappa_p(\blm,\bar \mu) = \kappa(\blm + \bar \mu) + \kappa(\blm) + \kappa(\bar\mu)$
is biadditive.  In that case, $\kappa_p$ is called the \emph{polarization}
of $\kappa$.  (The usual assumptions that
$\kappa(a\blm) = a^2\kappa(\blm)$ for $a\in \Ztwo$ and that $\kappa_p$ is $\Ztwo$-bilinear are automatic.)
If $\Lm$ has basis $\set{\lm_1,\dots,\lm_n}$, then the quadratic forms on
$\bLm$ are precisely the maps $\kappa : \bLm  \to \Ztwo$ of the form
\begin{equation}
\label{eq:Crquadratic}
\kappa(\sum_{i=1}^n l_i\blm_i) = \sum_{i=1}^n {l_ib_i }
+ \sum_{1\le i \le j \le n} l_il_ja_{ij},
\end{equation}
where $(b_i)\in \Ztwo^n$ and $(a_{ij})$ is an alternating matrix
in $\Mat_n(\Ztwo)$ (that is $a_{ij} = a_{ji}$ and $a_{ii} = 0$).

\begin{proposition}
\label{prop:Crquadratic}
(i)  If
$(\rA,\iota)$ is an associative $\Lm$-torus with involution,
then  there is a unique quadratic
form $\kappa : \bLm \to \Ztwo$, which we call the
\emph{mod-2 quadratic form for $(\rA,\iota)$}\index{associative torus with involution!mod-2 quadratic form},
such that
\begin{equation}
\label{eq:definekappa}
\iota(x_\lm) = (-1)^{\kappa(\blm)} x_\lm
\end{equation}
for $x_\lm\in \rA^\lm$, $\lm\in\Lm$.  In that case, we also have
\begin{equation}
\label{eq:propkappa}
x_\lm x_\mu = (-1)^{\kappa_p(\blm,\bmu)} x_\mu x_\lm
\end{equation}
for $x_\lm\in\rA^\lm$, $x_\mu\in\rA^\mu$, $\lm,\mu\in \Lm$.

(ii) Any quadratic form $\kappa$ on $\bLm$ is the mod-2 quadratic form for some
associative $\Lm$-torus with involution. In fact,  if $\Lm = \Zn$ with standard
basis $\set{\lm_1,\dots,\lm_n}$, $\kappa$
is given by \eqref{eq:Crquadratic} and we define $\be$ and $\bq$ by
$e_i = (-1)^{b_i}$ and $q_{ij} = (-1)^{a_{ij}}$, then $\kappa$
is the mod-2 quadratic form for $(\base_\bq,\iota_\be)$.

(iii) If $(\rA,\iota)$ is an associative $\Lm$-torus with involution and
$(\rA',\iota')$ is an associative $\Lm'$-torus with involution, then
$(\rA,\iota)\ig (\rA',\iota')$ if and only the corresponding
mod-2 quadratic forms $\kappa$ and $\kappa'$ are  isometric, which means
that there is an $\Ztwo$-linear isomorphism $\bar \tau : \bLm \to \bLm'$
such that $\kappa'(\bar\tau(\blm)) =  \kappa(\blm)$ for $\blm\in \bLm$.
\end{proposition}

\begin{proof}  We may assume that $\Lm = \Zn$ with standard basis $\set{\lm_1,\dots,\lm_n}$,
and in (iii) that $\Lm' = \Zn$.

(i) The uniqueness of $\kappa$ is clear.
For the existence we may assume that $(\rA,\iota) = (\base_\bq,\iota_\be)$
as in Example \ref{ex:quantuminvolution}.  Choose
$(b_i) \in \Ztwo^n$ and alternating $(a_{ij})\in\Mat_n(\Ztwo)$ such that
$e_i = (-1)^{b_i}$ and
$q_{ij} = (-1)^{a_{ij}}$. Then an easy computation
shows that  \eqref{eq:definekappa} holds with $\kappa$ given by
\eqref{eq:Crquadratic}.
Also, if  $x_\lm\in\rA^\lm$, $x_\mu\in\rA^\mu$, $\lm,\mu\in \Lm$,
then
$x_\lm x_\mu =  \iota( \iota(x_\mu)\iota(x_\lm))
= (-1)^{\kappa(\blm + \bmu)+ \kappa(\blm)+ \kappa(\bmu)}x_\mu x_\lm
$, so \eqref{eq:propkappa} holds.

(ii) It is enough to prove the second statement which follows from the proof of~(i).

(iii) We may assume that $(\rA,\iota) = (\base_\bq,\iota_\be)$
as in Example \ref{ex:quantuminvolution}.
For ``$\Rightarrow$'', let
$\eta: (\rA,\iota) \to (\rA',\iota')$ be an isograded-isomorphism
and set $\tau := \eta_\grtwist\in \GL(\Lm)$.  Then, $\eta(\rA^\lm) = {\rA'}^{\tau(\lm)}$
and so $\kappa'(\overline{\tau(\lm)}) = \kappa(\blm)$ for $\lm\in\Lm$.
Thus, the map $\bar \tau \in \GL(\bLm)$ induced by $\tau$ is an isometry of
$\kappa$ onto $\kappa'$.  For ``$\Leftarrow$'', let $\bar\tau\in \GL(\bLm)$
be an isometry of $\kappa$ onto $\kappa'$.  Then $\bar\tau$ is induced
by some $\tau\in \GL(\Lm)$.   (This follows from the fact that any
element of $\GL_n(\Ztwo)$ is the product of elementary matrices.)
Let $\lm_i' = \tau(\lm_i)$ and  $0\ne x_i' \in {\rA'}^{\lm_i'}$ for each $i$.
Then, $x_j'x_i' = (-1)^{\kappa_p'(\overline{\lm_i'},\overline{\lm_j'})}x_i'x_j'
= (-1)^{\kappa_p(\blm_i,\blm_j)}x_i'x_j'= q_{ij}x_i'x_j'$
and $\iota'(x_i') = (-1)^{\kappa'(\overline{\lm_i'})}x_i'
= (-1)^{\kappa(\blm_i)}x_i' = e_i x_i'$, so
$(\rA',\iota') \ig  (\base_\bq,\iota_\be)$ (see Example~\ref{ex:quantuminvolution}).
\qed\end{proof}

Finite dimensional quadratic forms over $\Ztwo$ have
been classified up to isometry \cite[Chap.~I, \S 16]{D},
and hence Proposition \ref{prop:Crquadratic} gives a corresponding
classification of  associative $\Lm$-tori with involution.

\begin{remark}
\label{rem:Crsupport}
Suppose that $(\rA,\iota)$ is an associative
$\Lm$-torus with involution  with $\rA_+ =  \rA_{+,\,\iota} $,
and set $\Lm_+ = \Lm_+(\rA,\iota) := \set{\lm\in\Lm \suchthat \rA^\lm \subseteq \rA_+}$.
If $1\le i \ne j \le r$, then
$\ell_{ij}(\rA)$, $m_{ij}(\rA)$, and $n_{ij}(\rA)$ are subspaces of $\rL(\rA,\iota)$
and the maps
$x\mapsto \ell_{ij}(x)$, $x\mapsto m_{ij}(x)$ and $x\mapsto n_{ij}(x)$ are injections from
$\rA$ onto $\ell_{ij}(\rA)$, $m_{ij}(\rA)$, and $n_{ij}(\rA)$ respectively.
On the other hand
if $1\le i \le r$, we have $m_{ii}(x) =  m_{ii}(\bar x)$ and
$n_{ii}(x) =  n_{ii}(\bar x)$ in $\rL(\rA,\iota)$ for $x\in \rA$,  and the maps
$h \mapsto m_{ii}(h)$ and $h \mapsto n_{ii}(h)$ are injections of
$\rA_+$ onto $m_{ii}(\rA_+) = m_{ii}(\rA)$ and $n_{ii}(\rA_+)= n_{ii}(\rA)$ respectively.
Thus, if $\al\in\Dec$, we have
\[\Lm_\al(\rL(\rA,\iota)) = \left\{
  \begin{array}{ll}
    \Lm, & \hbox{if $\al$ is short;} \\
    \Lm_+, & \hbox{if $\al$ is long.}
  \end{array}
\right.\]
Hence, if $\shom\in\Hom(Q,\Lm)$, then $\shom$ is admissible for $\rL(\rA,\iota)$
if and only if $\shom(\al_r)\in \Lm_+$.
\end{remark}

If $(\rA,\iota)$ is
an associative $\Lm$-torus with involution,
elements of $\rA_{+,\,\iota}$ are called \emph{hermitian} elements of $(\rA,\iota)$.
To form an isotope of $(\rA,\iota)$ along the lines of previous cases, we could a
define new product $(x,y)\mapsto xhy$ on $\rA$, where $h$ is a nonzero homogeneous
hermitian element of $(\rA,\iota)$, shift the grading, but keep the same involution.
Following \cite[\S I.3.4]{Mc2}, an equivalent but simpler approach is
to keep the product and grading and modify the  involution.

\begin{definition}[Isotopes of involutions]
\label{def:Crisotope}
Suppose that ($\rA,\iota)$ is
an associative $\Lm$-torus with involution.   Let
$h$ be a nonzero homogeneous hermitian element of $(\rA,\iota)$, and define
$\iota^\hi : \rA \to \rA$ by
\[\iota^\hi (x) = h\bar x h^{-1}\]
for $x\in \rA$.  Then,
$(\rA,\iota^\hi)$\index{associative torus with involution!isotope!$(\rA,\iota^\hi)$}
is an
associative $\Lm$-torus with involution.  (The multiplication
and grading on $\rA$ are unchanged.)  We call
$\iota^\hi$ the \emph{$h$-isotope} of~$\iota$.
\index{associative torus with involution!isotope}
\end{definition}

\begin{lemma}[Key lemma for type C$_r$]
\label{lem:Crkey}
Suppose $(\rA',\iota')$ is an associative $\Lm'$-torus with involution, $(\rA,\iota)$ is
an  associative $\Lm$-torus with involution, and $\shom\in \Hom(Q,\Lm)$ is admissible for
$\rL(\rA,\iota)$.  Suppose that $\ph : \rL(\rA',\iota') \to \rL(\rA,\iota)^\si$ is a bi-isomorphism
with $\ph_r =1$.  Choose $0\ne h\in \rA^{\shom(\al_r)}$.
Then there is an isograded-isomorphism $\eta : (\rA',\iota') \to (\rA,\iota^\hi)$
such that $\eta_\grtwist = \ph_e\in \GL(\Lm)$.
\end{lemma}

\begin{proof} In the argument $i,j,k$ will always denote integers  with $1\le i,j,k \le r$.

As in the proof of Lemma \ref{lem:A2key}, we see that
$\ph(\ell'_{ij}(\rA'^\lm)) = \ell_{ij}(\rA^{\ph_e(\lm) + \shom(\ep_i-\ep_j)})$
and
$\ph(m'_{ij}(\rA'^\lm)) = m_{ij}(\rA^{\ph_e(\lm) + \shom(\ep_i+\ep_j)})$
for $\lm\in\Lm$ and $i\ne j$.
So we may define linear bijections $\ph_{ij} : \rA'\to \rA$ and $\psi_{ij} : \rA' \to \rA$
for $i\ne j$ such that
\[
\ph(\ell_{ij}(x)) = \ell_{ij}(\ph_{ij}(x)),\quad \ph(m_{ij}(x)) = m_{ij}(\psi_{ij}(x))
\]
for $x\in \rA'$, and
\[
\ph_{ij}({\rA'}^\lm) =  \rA^{\ph_e(\lm)+ \shom(\ep_i-\ep_j)}, \quad
\psi_{ij}({\rA'}^\lm) =  \rA^{\ph_e(\lm)+ \shom(\ep_i+\ep_j)}
\]
for $\lm\in\Lm$.

Next, direct multiplication in $\rL(\rA,\iota)$ gives the following identities:
\[[\ell_{ik}(x),\ell_{kj}(y)] = \ell_{ij}(xy) \andd
[\ell_{ij}(x),[\ell_{ji}(y),m_{ij}(z)]] = m_{ij}(x(yz + \overline{yz}))
\]
for $x,y,z\in\rA$, $i,j,k\ne$.
We have similar equations in $\rL(\rA',\iota')$, and applying $\ph$ to these yields
\begin{align}
\ph_{ij}(xy) &= \ph_{ik}(x)\ph_{kj}(y),  \label{eq:Cr1}\\
\psi_{ij}(x(yz + \overline{yz})) &=
\ph_{ij}(x)(\ph_{ji}(y)\psi_{ij}(z)+\overline{\ph_{ji}(y)\psi_{ij}(z)}\ )
\label{eq:Cr2}
\end{align}
for $x,y,z\in\rA'$, $i,j,k\ne$.

We let
$u_{ij} = \ph_{ij}(1')\in \rA^{\shom(\ep_i-\ep_j)}$ and
$v_{ij} = \psi_{ij}(1') \in \rA^{\shom(\ep_i+\ep_j)}$ for $i\ne j$.
Setting $x=y=1'$ in \eqref{eq:Cr1}, we get $u_{ij} = u_{ik}u_{kj}$  for $i,j,k\ne$.
Therefore, $u_{ij}u_{ji} = u_{ki}^{-1}u_{ki}u_{ij}u_{ji} = u_{ki}^{-1}u_{kj}u_{ji} = u_{ki}^{-1}u_{ki} = 1$,
so $u_{ij} = u_{ji}^{-1}$ for $i\ne j$.

As in the proof of Lemma \ref{lem:A2key}, it follows from \eqref{eq:Cr1} that  we have
$\ph_{ij}(xy) = (\ph_{ij}(x)u_{jk})(u_{ki}\ph_{ij}(y))$ for $x,y\in \rA'$, $i,j,k\ne$.
So $\ph_{ij}(xy) = \ph_{ij}(x)u_{ji}\ph_{ij}(y)$ for $i\ne j$.
Thus, we see that $\chi_{ij}: \rA' \to \rA$, defined by
\begin{equation}
\label{eq:chiijdef}
\chi_{ij}(x) = u_{ji}\ph_{ij}(x),
\end{equation}
is an algebra isomorphism of $\rA'$ onto $\rA$ for $i\ne j$.
Also, $\chi_{ij}(\rA'^\lm) =  u_{ji}\ph_{ji}(\rA'^\lm) =
u_{ji} \rA^{\ph_e(\lm)+ \shom(\ep_i-\ep_i)} \subseteq \rA^{\ph_e(\lm)}$
for $\lm\in \Lm$, so $\chi_{ij}$ is isograded with $(\chi_{ij})_\grtwist = \ph_e$.

If we put $x=y=1'$ in \eqref{eq:Cr2}, we get
$\psi_{ij}(z + \bar z) =
u_{ij}(u_{ji}\psi_{ij}(z)+\overline{u_{ji}\psi_{ij}(z)}\ )$.  Since
$u_{ij}u_{ji} = 1$, this yields
\begin{equation}
\label{eq:Cr3}
\psi_{ij}(\bar z) =
u_{ij}\, \overline{\psi_{ij}(z)} \, \bar u_{ji}
\end{equation}
for $z\in\rA'$, $i\ne j$.
Also, if we take $y=z=1'$  in \eqref{eq:Cr2}, we get
\begin{equation}
\label{eq:Cr4}
\psi_{ij}(2x) =
\ph_{ij}(x)(u_{ji}v_{ij} +\overline{u_{ji}v_{ij}})
\end{equation}
for $x\in \rA'$, $i\ne j$.

Next $u_{r1} v_{1,r}\in \rA^{\shom(\ep_r-\ep_1)} \rA^{\shom(\ep_r+\ep_1)} \subseteq
\rA^{\shom(2\ep_r)}$.  So replacing $h$ by a scalar multiple, we can assume that
$h = u_{r1} v_{1,r}$.  Since $\shom$ is admissible for $\rL(\rA,\iota)$,
we have $\bar h = h$, so \eqref{eq:Cr4} tells us that
\begin{equation}
\label{eq:Cr5}
\psi_{1r}(x) = \ph_{1r}(x)h
\end{equation}
for $x\in \rA'$.
Finally, using \eqref{eq:chiijdef}, \eqref{eq:Cr3} and \eqref{eq:Cr5}, we have
\begin{align*}
\chi_{1r}(\bar x) &= u_{r1}\ph_{1r}(\bar x) =
 u_{r1}\psi_{1r}(\bar x) h^{-1}
= u_{r1}u_{1r}\, \overline{\psi_{1r}(x)} \, \bar u_{r1} h^{-1}\\
&= \overline{\psi_{1r}(x)} \, \bar u_{r1} h^{-1}
= h \overline{\ph_{1r}(x)} \, \bar u_{r1} h^{-1}
= h \overline{\chi_{1r}(x)}\, h^{-1}
\end{align*}
for $x\in \rA'$.  Thus, $\chi_{1r}$ is an isograded  isomorphism of $(\rA',\iota')$
onto $(\rA,\iota^\hi)$.
\qed\end{proof}

As in types A$_1$ and A$_2$, the following proposition now follows from
Theorem \ref{thm:Crcoord} and Lemma \ref{lem:Crkey}.

\begin{proposition}
\label{prop:Crisotope}  Suppose that $(\rA,\iota)$ is
associative $\Lm$-torus with involution
and $\shom\in \Hom(Q,\Lm)$ is admissible for $\rL(\rA,\iota)$.
Choose $0\ne h\in \rA^{\shom(\al_r)}$. Then, we have
$ \rL(\rA,\iota)^\si \simeq_{Q\times\Lm}\rL(\rA,\iota^\hi)$.
\end{proposition}

We could follow the pattern in type A$_1$ and define two associative
tori with involution  $(\rA,\iota)$ and $(\rA',\iota')$ to be
isotopic\index{associative torus with involution!isotopic}
if
$(\rA,\iota^\hi)$ is isograded-isomorphic to $(\rA',\iota')$ for some
nonzero homogeneous hermitian element $h$ of $(\rA,\iota)$.  However,
we don't need to do this because of the following:

\begin{proposition}
\label{prop:Crisotopytrivial}
Suppose that $(\rA,\iota)$ is an associative  $\Lm$-torus with involution
and $h$ is a nonzero homogeneous hermitian
element of $(\rA,\iota)$. Then $(\rA,\iota^\hi) \ig (\rA,\iota)$.
\end{proposition}

\begin{proof} Let $\kappa$ and $\kappa^\hi$ be the mod-2 quadratic forms
for $(\rA,\iota)$ and $(\rA,\iota^\hi)$ respectively.  Now $h\in \rA^\mu$, where
$\mu\in \Lm$, and since $h$ is hermitian we have $\kappa(\bar \mu) = 0$.
Also, if $x\in \rA^\lm$, $\lm\in \Lm$, we have
$\iota^\hi(x) = h \iota(x)h^{-1} = (-1)^{\kappa(\blm)} h x h^{-1}
= (-1)^{\kappa(\blm)+\kappa_p(\blm,\bmu)} x$, so
$\kappa^\hi(\blm) = \kappa(\blm)+\kappa_p(\blm,\bmu)$ for $\blm\in \bLm$.
Hence, since $\kappa(\bmu) = 0$, we have
\[\kappa^\hi(\blm) =  \kappa(\blm)+\kappa_p(\blm,\bmu)^2 = \kappa(\blm + \kappa_p(\bmu,\blm)\bmu)\]
for $\blm\in \bLm$.
So we define $\bar\tau : \bLm \to \bLm$ by $\bar\tau(\blm) = \blm + \kappa_p(\bmu,\blm)\bmu$
for $\blm\in \bLm$.  One checks that $\bar\tau^2 = 1$, so $\bar\tau$ is an isometry
of $\kappa^\hi$ onto $\kappa$.  Hence, our conclusion follows from
Proposition \ref{prop:Crquadratic}(iii).
\qed\end{proof}

\begin{theorem}
\label{thm:Crisotopy}
Suppose that $(\rA,\iota)$ is an associative $\Lm$-torus with involution
and $(\rA',\iota')$ is an associative $\Lm'$-torus with involution.
Let $\rL(\rA,\iota) = \ssp_{2r}(\rA,\iota)$ and $\rL(\rA',\iota') =  \ssp_{2r}(\rA',\iota')$.
Then,
the following are equivalent:
\begin{description}[(a) ]
\item[(a)] $(\rA,\iota)$ is isograded-isomorphic to $(\rA',\iota')$.
\item[(b)] $\rL(\rA,\iota) \bi \rL(\rA',\iota')$.
\item[(c)] $\rL(\rA,\iota) \sim \rL(\rA',\iota')$.
\end{description}
\end{theorem}

\begin{proof}  ``(a)$\Rightarrow$(b)'' is clear, and
``(b)$\Rightarrow$(c)'' is trivial. Finally, suppose that
(c) holds. Thus,
$\rL(\rA',\iota') \bi \rL(\rA,\iota)$ by Propositions \ref{prop:Crisotope} and  \ref{prop:Crisotopytrivial},
so there exists a bi-isomorphism $\ph : \rL(\rA',\iota') \to \rL(\rA,\iota)$.
Since the automorphism group of $\De$ equals the Weyl group of $\De$,
we can, by Lemma \ref{lem:bi-isomorphism}, assume that $\ph_r = 1$.
(a) now follows from Proposition \ref{prop:Crisotope} (with $\shom = 0$
and  $h = 1$).
\qed\end{proof}

\section{Some concluding remarks}
\label{sec:conclusions}

Putting together the theorems quoted and proved in this article, we obtain a very good understanding
of the structure of EALAs of several types.

\medskip\par\noindent\emph{Type A$_1$:} By Theorems \ref{thm:A1coord} and Theorem \ref{thm:A1isotopy},
the TKK construction induces a 1-1 correspondence between isotopy classes
of Jordan $\Zn$-tori and isotopy classes of centreless Lie $\Zn$-tori of type A$_1$.
Composing this with Neher's EALA construction (Construction \ref{con:famEALA}),
we obtain, by the results of \S  \ref{sec:isotopyEALA},
an explicit 1-1 correspondence between isotopy classes of
Jordan $\Zn$-tori and families of EALAs of type A$_1$ and nullity $n$
up to bijective isomorphism.

\medskip\par\noindent\emph{Type A$_r\,  (r \ge 2)$:} By the results of
\S's \ref{sec:isotopyEALA}, \ref{sec:A2} and \ref{sec:Ar},
if $r = 2$ (resp.~$r\ge 3$) the $3\times 3$-projective elementary construction (resp.~the $(r+1)\times (r+1)$-special linear construction) followed by Neher's construction gives us an explicit 1-1 correspondence between isograded-isomorphism classes of alternative (resp.~associative) $\Zn$-tori and families of EALAs of type A$_r$
and nullity $n$ up to bijective isomorphism.

\medskip\par\noindent\emph{Type C$_r\,  (r \ge 4)$:} The structural information is particularly sharp in this case.
Indeed, by Theorem \ref{thm:Crcoord} and Theorem \ref{thm:Crisotopy}
(again along with the results of \S \ref{sec:isotopyEALA}), the
construction of an  associative torus with involution from a mod-2 quadratic form (see Proposition
\ref{prop:Crquadratic}(ii)),
followed by the $2r\times 2r$-special symplectic construction,
followed by Neher's construction, gives an explicit 1-1 correspondence between
isometry classes of quadratic forms on an $n$-dimensional vector space over $\Ztwo$ and
families of EALAs of type C$_r$ and nullity $n$ up to bijective isomorphism.

\medskip\par\noindent\emph{Type D$_r\, (r \ge 4)$, E$_6$, E$_7$, E$_8$:}  Let X$_r$
be one of the indicated types.  According to
\cite[Thm.~1.37]{BGK} and \cite[\S 5]{Y2},  any centreless Lie $\Lm$-torus of type X$_r$ is
isograded-isomorphic to the untwisted Lie torus $\fg\otimes \base[\Lm]$, where
$\fg$ is the split simple Lie algebra of type X$_r$. Thus, there is only
one bi-isomorphism class and hence only one isotopy class
of centreless Lie $\Zn$-tori of type X$_r$.  Therefore, for a given nullity,
there is only one family of EALAs of type X$_r$ up to bijective isomorphism.

\medskip\par\noindent\emph{Other types:}
We have omitted discussion of families of EALAs of types C$_3$,
B$_r\, (r\ge 2)$  and BC$_r\, (r\ge 1)$.  Type C$_3$ can be handled
using for the methods of \S \ref{sec:Cr},
as there is only one centreless Lie $\Zn$-torus of type C$_3$ that
is not covered by the $6\times 6$-special symplectic construction
(\cite[Thm.~4.87]{AG}, \cite[Thm.~5.9]{BY}).
(This exceptional example exists only if $n\ge 3$ .)
In types B$_r\, (r\ge 2)$  and BC$_r\, (r\ge 1)$, we expect that isotopes
and isotopy of coordinate algebras will play a key role not only in the arguments but also
in the statements of the results (as in type A$_1$).

\medskip
Finally, we note that when $\base$ is algebraically closed, there is an
alternative method (replacing coordinatization) to construct centreless Lie $\Zn$-tori.
This me\-thod, which has been investigated recently in \cite{ABFP1} and \cite{ABFP2},
constructs centreless Lie tori as multiloop algebras, just (as in nullity 1) affine Kac-Moody
Lie algebras are constructed using loop algebras. To be somewhat more precise,
suppose that $\bsg = (\sg_1,\dots,\sg_n)$ is an $n$-tuple of commuting finite
order automorphisms of a finite dimensional simple Lie algebra $\fs$.
Let $\fg = \fs^\bsg$ be the fixed point algebra of $\set{\sg_i}_{i=1}^n$  in $\fs$,
let $\fh$ be a Cartan subalgebra of $\fg$, let
$\De = \De(\fs,\fh)$, and let $Q = \spann_\bbZ{\De}$.
We give the tensor product $\rT = \fs\ot \base[z_1^{\pm 1},\dots,z_n^{\pm 1}]$
a $Q\times \Zn$-grading using the natural $Q$-grading
of the first factor and the natural  $\Zn$-grading of the second factor.
Next, for $1\le i \le n$,
we fix a primitive $m_i^\text{th}$ root of unity
$\zeta_{m_i}$ for each $i$, where $m_i = \order{\sg_i}$, and we define
an automorphism $\eta_i$ of
$\base[z_1^{\pm 1},\dots,z_n^{\pm 1}]$
by $\eta_i(z_j) = \zeta_{m_i}^{\delta_{ij}}z_j$.
Now let $\LT(\fs,\bsg,\fh)$ be the
fixed point algebra of $\set{\sg_i \ot \eta_i^{-1}}_{i=1}^n$  in $\rT$,
in which case $\LT(\fs,\bsg,\fh)$ is a $Q\times \Zn$-graded subalgebra of $\rT$.
If $\bsg$ satisfies some additional conditions that are listed in \cite[Prop.~3.2.5]{ABFP2}
(the most important being that $\fg$ is simple),
then $\De$ is a root system in $\fh^*$ and $\LT(\fs,\bsg,\fh)$ is a centreless Lie
$\Zn$-torus of type $\De$, which is called a
\emph{multiloop Lie torus}\index{Lie torus!multiloop}.
Also it is shown in \cite[Thm.~3.3.1]{ABFP2} that
a centreless Lie torus $\rL$ is bi-isomorphic to some multiloop Lie torus if and only if
$\rL$ is a finitely generated as a module over its centroid.  Furthermore,
necessary and sufficient conditions are given there for two multiloop Lie tori
to be bi-isomorphic or to be isotopic.
The reader can consult \cite{ABFP1} and \cite{ABFP2} for details about these results.

The two approaches, using coordinatization and using multiloops, which have been used
to construct
and to determine isotopy of centreless Lie tori,  each have significant
complimentary advantages.  It is our view that both are needed in
order to gain a full understanding of families of EALAs.

\printindex

\end{document}